\documentclass[11pt]{article}
\usepackage{graphicx} 

\usepackage[letterpaper, margin=1in]{geometry}
\usepackage{amsthm,amsmath}
\usepackage{natbib,setspace}
\usepackage{color,comment}

\newtheorem{theorem}{Theorem}
\newtheorem{corollary}{Corollary}
\newtheorem{lemma}{Lemma}
\newtheorem{assumption}{Assumption}

\usepackage{amsmath,amssymb,amsfonts}

\usepackage{natbib}
 \bibpunct[, ]{(}{)}{,}{a}{}{,}%

\usepackage{rotating}
\usepackage{fancyvrb}
\usepackage{mathrsfs}
\newcommand{\ve}{\mathscr{V}}

\begin{document}



\title{Rates of Convergence in the Central Limit Theorem for Markov
Chains, with an Application to TD Learning}


\author{R. Srikant\\ECE and CSL, UIUC\\{rsrikant@illinois.edu}}


\maketitle

\begin{abstract}
     We prove a non-asymptotic central limit theorem for vector-valued martingale differences using Stein's method, and use Poisson's equation to extend the result to functions of Markov Chains. We then show that these results can be applied to establish a non-asymptotic central limit theorem for Temporal Difference (TD) learning with averaging.
\end{abstract}

\section{Introduction}

Starting with the seminal work of \citep{polyak1992acceleration}, where they analyzed the averaging rule of \citep{polyak1990new,ruppert1988efficient}, the central limit theorem for martingales has been a valuable tool to understand asymptotic efficiency of stochastic approximation algorithms for optimization. More recently, the central limit theorem for functions of Markov chains has also been used to study the asymptotic efficiency of machine learning algorithms \citep{borkar2021ode,hu2024central}. However, in machine learning applications, finite-time (rather than asymptotic) bounds are preferred to understand the sample complexity of algorithms. However, to the best of our knowledge, there is limited work on the central limit theorems for vector-valued martingales and vector-valued functions of Markov chains that would allow one to both quantify asymptotic efficiency and obtain finite-time bounds for learning and optimization algorithms. A recent paper by \citep{anastasiou2019normal} derives a rate of convergence for the central limit theorem vector-valued martingales; however, the notion of distance between probability distributions in that paper is not sufficient for some applications.

Our goal in this paper is to obtain bounds on the rate of convergence to normality to complement the central limit theorem for Markov chains \citep{meyn2012markov}. To this end, we first establish a result for martingale differences by combining the Lindeberg decomposition along with Stein's method, which was first used in \citep{rollin2018quantitative} for scalar valued martingale differences. This key idea in \citep{rollin2018quantitative} was extended to the case of vector-valued martingales in \citep{anastasiou2019normal} who also used it to obtain a non-asymptotic version of the central limit theorem in \citep{polyak1992acceleration}. Our approach builds upon the ideas in \citep{rollin2018quantitative,anastasiou2019normal}, but we also introduce new techniques in order to obtain stronger results:
\begin{itemize}
    \item We would like our results to hold when the deviation from normality is measured using the Wasserstein distance. The results in \citep{anastasiou2019normal} provide bounds only for weaker distances. Therefore, we have to obtain results for more general test functions than the ones considered in \citep{anastasiou2019normal}. To do so, we leverage the results on the regularity of solutions to Stein's equation for vector-valued random variables in \citep{gallouet2018regularity,fang2019multivariate}, which makes our analysis considerably different. 
    \item To obtain rates of convergence for the Markov chain central limit theorem, we convert a Markov chain into a martingale difference noise using Poisson's equation; see \citep{metivier1984applications,benveniste2012adaptive,douc2018markov}, for example. To the best of our knowledge, this idea has not been used to derive convergence rates for the central limit theorem for functions of Markov chains. 
\end{itemize}

A difference between our results and the results in \citep{rollin2018quantitative,anastasiou2019normal} for martingales is the fact that we state our convergence bounds in terms of a positive definite matrix $\Sigma_\infty$ which can be viewed as the limiting expected conditional covariance of the martingale differences (see Theorem~\ref{thm: martingale CLT} and the short discussion after the theorem). This is important to us since our goal is to state our results for Markov chains in terms of the asymptotic covariance, which is important in applications; see \citep{borkar2021ode}, for example. Therefore, our assumptions are stated differently than in \cite[Theorem 1]{rollin2018quantitative} or in \cite[Equation 3.1]{anastasiou2019normal}. It is also worth remarking that Stein's method to derive rates of convergence for the martingale central limit theorem is relatively new. The early papers on obtaining rates of convergence for the martingale central limit theorem use other techniques; see, for example, \citep{bolthausen1982exact} and other references in \citep{rollin2018quantitative}. Here, we use Stein's method because we found it easier to quantify the impact of the asymptotic variance on the rate of convergence using Stein's method. Specifically, we show that the regularity results in \citep{gallouet2018regularity,fang2019multivariate} and the calculations in \citep{rollin2018quantitative,gallouet2018regularity} can be used, along with insights from Markov chain theory, to understand the impact of the asymptotic covariance term on the rates of convergence in the martingale and Markov chain central limit theorems. 

As an application, we apply the Markov chain central limit theorem to study the distribution of error in TD learning. TD learning is the most widely studied algorithm for evaluating the performance of a given policy from data in Markov Decision Processes \citep{sutton2018reinforcement}. The convergence of TD learning with decaying stepsizes was proved in \citep{tsitsiklis1996analysis}. Over the last few years, there has been a resurgent interest in understanding the non-asymptotic convergence behavior of TD learning and further using it to study other reinforcement algorithms which use TD learning within their framework. The first result of this type was obtained in \citep{bhandari2018finite} who used a projection step in their version of the TD learning algorithm. The first finite-time bounds for unprojected TD learning were obtained in \citep{srikant2019finite}, who studied both fixed step-sizes and decaying step-sizes. TD learning can be viewed as a special case of linear stochastic approximation with Markovian and multiplicative noise and as a special case of Markovian jump linear systems \citep{hu2019characterizing}. The results in \citep{srikant2019finite} were extended to nonlinear, contractive stochastic approximation algorithms in \citep{chen2023lyapunov}, who also derived finite-time performance bounds for specific choices of decaying step-size rules, which also allowed them to study the robustness of different step-size choices. Finite-time performance bounds for variants of TD learning and other schemes have been studied in \citep{gupta2019finite,doan2021finite,doan2022nonlinear} and concentration results for TD learning have been obtained in \citep{telgarsky2022stochastic,chen2023concentration}.

While most of the above papers study standard TD learning, in practice, it is well-known that averaging the iterates can help improve the asymptotic variance. This idea was originally studied for stochastic approximation and stochastic gradient descent in the seminal works in \citep{ruppert1988efficient,polyak1990new,polyak1992acceleration}; see \citep{moulines2011non} for a finite-time analysis. More recently, the asymptotic behavior of TD learning and other algorithms with averaging has been studied in \citep{borkar2021ode,mou2021optimal,durmus2022finite}; see also \citep{li2023online} for related work. As in the papers on stochastic approximation, the conclusion in \citep{borkar2021ode} is that, for a suitable choice of decaying stepsizes, averaging the output leads to efficient estimators in the sense that the variance is minimized. Further, \citep{borkar2021ode} also established functional central limit theorems and central limit theorems for TD learning and other reinforcement learning algorithms. The work in \citep{mou2021optimal} shows order-optimal rate of convergence for the mean-squared error, while \citep{durmus2022finite} obtains higher moment bounds and high probability deviation bounds. In contemporaneous work to this paper, finite-time performance bounds on the covariance of the averaged TD learning iterates have been obtained in \citep{haque2023tight} using more robust step-size rules. Their model is more general; they study two-time scale linear stochastic approximation algorithms, with TD learning with averaging being a special case of the models that they study. A central limit theorem for two time-scale stochastic approximation has also been established for two time-scale stochastic approximation in \citep{hu2024central}. Related work also includes a central limit theorem for Polyak-Ruppert averaged Q-learning \citep{li2023statistical}.

Our main results in the paper are as follows\footnote{This version of the paper corrects some errors in the version published in Math of OR; while the list of contributions is the same, updated theorems and proofs are presented later in the paper.}:
\begin{itemize}
    \item We derive non-asymptotic central limit theorems for martingales and Markov chains by exploiting ideas from \citep{rollin2018quantitative,gallouet2018regularity,fang2019multivariate}. 
    \item We then extend the results to Markov chains using Poisson's equation to relate Markov chains to martingales. In this case, the asymptotic covariance can be characterized as is well known in the literature; see \citep{douc2018markov}, for example. 
    \item Finally, we apply the Markov chain results, along with the ideas in \citep{polyak1992acceleration}, to estimate the rate of convergence in the central limit theorem for TD learning with Polyak-Ruppert averaging. An advantage of the approach here is that our result directly establishes a central limit theorem for the scaled averaged of the iterates, instead of first establishing a function central limit theorem for the entire trajectory of the iterates and then establishing the central limit theorem.
    \item In addition to the difference in the notion of distance considered here and in \citep{anastasiou2019normal}, in the TD learning application, we have to deal with the fact that the noise multiplies the TD learning parameters, whereas the noise is additive in the stochastic gradient application considered in \citep{anastasiou2019normal}. 
\end{itemize}

It is important to note that, unlike \citep{polyak1992acceleration}, our result only holds for averaging along with an appropriate decaying step-size rule. Specifically, the result does not hold for fixed step-size TD learning with averaging due to the bias issue pointed out in \citep{borkar2021ode,lauand2023curse,huo2023bias}); also see \citep{nagaraj2020least,roy2023online,jain2018parallelizing} for challenges in dealing with Markov noise in other problems. A nonasymptotic central limit theorem has been established for constant step-size stochastic gradient descent in \citep{dieuleveut2020bridging}, but the model considered here is different from the problem in \citep{dieuleveut2020bridging} (for example, Markovian noise, constant step-size vs decaying step-size) and therefore, the techniques there do not seem to apply to our problem.

The rest of the paper is organized as follows. In Section~\ref{sec: Martingales}, we establish an upper bound on the rate of convergence in the martingale central limit theorem under some assumptions, and in Section~\ref{sec: Markov}, we apply these results to bound the rate of convergence in the central limit theorem for functions of Markov chains using Poisson's equation. The Markov chain results are then applied to study the central limit for averaged iterates of TD learning in Section~\ref{sec: TD}. Concluding remarks and suggestions for future work are provided in Section~\ref{sec: conclusions}.

\section{Martingale Central Limit Theorem}\label{sec: Martingales}

In this section, we present a bound on the rate of convergence for the martingale central limit theorem in terms of the Wasserstein distance using Stein's method.

\subsection{Preliminaries}

In this subsection, we will describe the notation we use to establish the rate of convergence in the central limit theorem for vector-valued martingale differences. Let $\{m_k\}_{k\geq 1}$ be a $d$-dimensional martingale difference sequence with respect to a filtration $\{\mathcal{F}_k\}_{k\geq 0}$ i.e., $m_k$ is $\mathcal{F}_k$-adapted, $E(||m_k||)<\infty,$ where $||\cdot||$ is the standard Euclidean norm in $\Re^d$ and $$E(m_k|\mathcal{F}_{k-1})=0\qquad \forall k\geq 1.$$ 
\begin{assumption}\label{assume: martingale}
Fix $\beta\in(0,1)$. We assume that the martingale difference sequence satisfies the following properties:
\begin{enumerate}\label{assume}
\item $E\|m_k\|^{2+\beta}<\infty$ for every $k\geq 1$.
\item $\Sigma_k=E(m_km_k^T\mid\mathcal{F}_{k-1})$ exists for every $k\geq 1$.
\end{enumerate}
\hfill $\diamond$
\end{assumption}
Item (1) in the assumption is motivated by the fact that we want to extend the techniques in \citep{gallouet2018regularity} to martingales .

Define
$S_n=\sum_{k=1}^n m_k.$ We are interested in whether the following normalized sum 
\begin{equation}\label{eq: normalized sum}
W_n:=\frac{S_n}{\sqrt{n}}
\end{equation}
converges to $N(0,\Sigma_\infty)$, and if so, we are interested in the rate of convergence. We will use Stein's method to establish such a convergence and to estimate the rate of convergence \citep{stein1972bound,ross2011fundamentals,chatterjee2014short}. 

The Wasserstein distance between two $d$-dimensional random vectors $X$ and $Y$ is defined as
$$d_{\mathcal{W}}(X,Y)=\sup_{h\in Lip_1} E[h(X)-h(Y)],$$
where $Lip_L$, the class of Lipschitz functions from $\Re^d$ to $\Re$  with Lipschitz constant $L,$ is defined as
$$Lip_L=\{h: |h(x)-h(y)|\leq L||x-y||\}.$$
For later reference, we use the notation $||A||_{op}$  to denote the operator norm of a matrix $A,$ i.e., $$||A||_{op}=\sup_{x: ||x||=1}||Ax||.$$
The following results from \cite[Propositions 2.2 and 2.3]{gallouet2018regularity} and \cite[Theorem 3.1]{fang2019multivariate} will be useful to us and so we present them together as a lemma for future reference.
\begin{lemma}\label{lem: regularity}
Let $g\in Lip_L$ and $f$ be the solution to Stein's equation:
\begin{equation}\label{eq: Stein}
    g(u)-E(g(Z))=\Delta f(u)-u^T\nabla f(u),
\end{equation}
where $\Delta$ and $\nabla$ are the Laplacian and gradient operators, respectively, and $Z\sim N(0,I)$. Then, $f$ has the following regularity properties for any $\beta\in (0,1)$:
\begin{enumerate}
    \item $||\nabla^2f(x)-\nabla^2f(y)||_{op}\leq \tilde{C}_1(d,\beta)||x-y||^\beta L,$ where $\nabla^2$ is the Hessian operator.
    \item $|\Delta f(x)-\Delta f(y)|\leq \tilde{C_2}(d,\beta)||x-y||^\beta L$,
where  $$\tilde{C}_1=C_1(d)+\frac{2}{1-\beta}, \quad \tilde{C}_2=C_2(d)+\frac{2d}{1-\beta}, \quad 
C_1=2^{3/2}\frac{(1+2d)\Gamma((1+d)/2)}{d\Gamma(d/2)},\quad C_2=2\sqrt{\frac{2d}{\pi}}.$$
    \item  $$\left|\langle\nabla^2 f(x),u_1u_2^T\rangle_{HS}\right|\leq C L||u_1||\cdot||u_2||\,\forall x, u_1, u_2\in \Re^d,$$ where $C$ is a constant (independent of $d$ and $\beta$) and $\langle\cdot,\cdot\rangle_{HS}$ is the Hilbert-Schmidt inner product (also called the Frobenius inner product), i.e., the sum of the element-by-element product of two matrices.
\end{enumerate}
\hfill $\diamond$
\end{lemma}
Next, we state a well-known result regarding Ornstein-Uhlenbeck processes which will be useful later to establish the nonasymptotic version of the martingale central limit theorem  (for example, the result can be easily inferred from the discussions on generators and linear stochastic differential equations in \citep{arnold1974stochastic}).
\begin{lemma}\label{lem: O-U}
    Consider the O-U process
    $$dX_t=AX_tdt+Bdw_t,$$
    where $w$ is the standard $d$-dimensional Wiener process, and $A$ and $B$ are $d\times d$ matrices. Then, the generator of the process $\mathcal{A}$ is the operator given by
    $$\mathcal{A}f(x):=\lim_{\delta\rightarrow 0}E\left(\frac{f(X_{t+\delta})-f(X_t)|X_t=x}{\delta}\right)=\nabla^Tf(x)Ax+\frac{1}{2}Tr(\nabla^2f(x) BB^T),$$
    where $f:\Re^d\rightarrow \Re$ is any twice differentiable function. If the eigenvalues of $A$ have strictly negative real parts, then the stationary distribution of the O-U process is Gaussian with mean zero and covariance $\Sigma,$ which satisfies the Lyapunov equation
    $$\Sigma A^T+A\Sigma+BB^T=0.$$ Thus, if $\tilde{Z}_k\sim N(0,\Sigma),$ then $E(\mathcal{A}f(\tilde{Z}_k))=0$ for any twice differentiable function $f.$
    \hfill $\diamond$
\end{lemma}
It is worth noting that Stein's equation (\ref{eq: Stein}) can be motivated by the above lemma when $A=-I$ and $BB^T=2I,$ in which case the right-hand side of (\ref{eq: Stein}) is the same as $\mathcal{A}f(u)$. Thus, if $u$ in (\ref{eq: Stein}) is replaced by a random vector, then the left-hand side of (\ref{eq: Stein}) is a measure of the difference between the distribution of that random vector and a standard Gaussian distribution while the right-hand side will hopefully be small if the distribution of the random vector is close to a standard Gaussian distribution (see, for example, \citep{barbour1990stein,gotze1991rate,chen2010normal}).

\subsection{Rate of Convergence}

We now state and prove the main result of this section.

\begin{theorem}\label{thm: martingale CLT}
Let
\[
    V_k:=\sum_{j=k}^n E(\Sigma_j),
    \qquad k=1,\ldots,n,
\]
and assume that \(V_k\) is positive definite for each \(k=1,\ldots,n\).  Define
\[
    D_k:=V_k^{-1/2}E(\Sigma_k)V_k^{-1/2},
    \qquad k=1,\ldots,n .
\]
Let \(Z\sim N(0,I)\), and let \(W_n\) be defined as in
\((\ref{eq: normalized sum})\). Then, under
Assumption~\ref{assume: martingale} for the same fixed \(\beta\in(0,1)\),
\begin{align*}
   d_{\mathcal W}(W_n,n^{-1/2}V_1^{1/2}Z)
   &\le
   \frac{1}{\sqrt n}\sum_{k=1}^n
   \Biggl[
   \tilde C_1(d,\beta)\|V_k^{1/2}\|_{op}
   E\left(\|V_k^{-1/2}m_k\|^{2+\beta}\right)
   \\
   &\qquad\qquad
   +
   \tilde C_1(d,\beta)\operatorname{Tr}(D_k)\|V_k^{1/2}\|_{op}
   E\left(\|V_k^{-1/2}m_k\|^\beta\right)
   \Biggr]
   +\Delta_n,
\end{align*}
where
\begin{align*}
\Delta_n
:=
\sup_{h\in Lip_1}
\frac{1}{\sqrt n}
\left|
\sum_{k=1}^n
E\left[
\operatorname{Tr}\left(
A_{k,h}(\mathcal F_{k-1})
V_k^{-1/2}
\left(E(\Sigma_k)-\Sigma_k\right)
V_k^{-1/2}
\right)
\right]
\right|.
\end{align*}
Here
\[
A_{k,h}(\mathcal F_{k-1})
:=
E\left[
\nabla^2 f_{k,h}(\widetilde Z_k)
\mid
\mathcal F_{k-1}
\right],
\]
where \(f_{k,h}\) is the $\mathcal F_{k-1}$-measurable random solution of
the Stein equation defined in the proof below with
\[
    a=S_{k-1},
    \qquad
    Q_k=V_k^{1/2}.
\]
and
\[
    \widetilde Z_k
    :=
    V_k^{-1/2}T_{k+1}
    \sim
    N\left(0,V_k^{-1/2}V_{k+1}V_k^{-1/2}\right),
\]
with the convention \(V_{n+1}=0\). Moreover,
\[
    \|A_{k,h}(\mathcal F_{k-1})\|_{op}
    \le
    C\|V_k^{1/2}\|_{op}, \quad \text{a.s.},
\]
where \(C\) is the constant in statement (3) of Lemma~\ref{lem: regularity}.
Since \(0\preceq D_k\preceq I\), the same bound also implies the coarser
version obtained by replacing \(\operatorname{Tr}(D_k)\) by \(d\).  The sharper
form above is used in the Markov-chain applications.
\end{theorem}

\begin{proof}
Fix \(h\in Lip_1\). By the standard rescaling of Lipschitz test functions,
\[
d_{\mathcal W}(W_n,n^{-1/2}V_1^{1/2}Z)
=
\sup_{h\in Lip_1}
\frac{1}{\sqrt n}
E\left[h(S_n)-h(T_1)\right],
\]
where
\[
    T_k:=\sum_{j=k}^n Y_j,
\]
and \(Y_1,\ldots,Y_n\) are independent Gaussian random vectors, independent
of the martingale difference sequence, with
\[
    Y_j\sim N(0,E(\Sigma_j)).
\]
Thus
\[
    T_k\sim N(0,V_k).
\]
The rescaling follows because if \(g\in Lip_1\), then
\(h(x):=\sqrt n\,g(x/\sqrt n)\) also belongs to \(Lip_1\), and
\[
E\left[g(W_n)-g(n^{-1/2}V_1^{1/2}Z)\right]
=
\frac{1}{\sqrt n}
E\left[h(S_n)-h(T_1)\right].
\]

Using the Lindeberg decomposition,
\[
E\left[h(S_n)-h(T_1)\right]
=
\sum_{k=1}^n
E\left[
h(S_k+T_{k+1})-h(S_{k-1}+T_k)
\right],
\]
where \(T_{n+1}=S_0=0\).

We now bound a single summand. Fix \(k\). Conditional on
\(\mathcal F_{k-1}\), set
\[
    a=S_{k-1},
    \qquad
    Q_k=V_k^{1/2}.
\]
Define
\[
    \widetilde h_k(u):=h(a+Q_k u).
\]
Since \(h\in Lip_1\), we have
\[
    \widetilde h_k\in Lip_{L_k},
    \qquad
    L_k=\|Q_k\|_{op}=\|V_k^{1/2}\|_{op}.
\]
Let \(f_k=f_{k,h}\) be the solution of the Stein equation
\[
    \widetilde h_k(u)-E[\widetilde h_k(Z)]
    =
    \Delta f_k(u)-u^\top\nabla f_k(u).
\]
Since \(a=S_{k-1}\) is \(\mathcal F_{k-1}\)-measurable, the function
\(f_k\) is also \(\mathcal F_{k-1}\)-measurable as a random function.

By the Stein equation, with
\[
    \widetilde Z_k:=V_k^{-1/2}T_{k+1},
\]
we have
\begin{align*}
& E\left[
h(S_k+T_{k+1})-h(S_{k-1}+T_k)
\mid
\mathcal F_{k-1}
\right]
\\
&=
E\biggl[
\Delta f_k\left(V_k^{-1/2}m_k+\widetilde Z_k\right)
-
\left(V_k^{-1/2}m_k+\widetilde Z_k\right)^\top
\nabla f_k\left(V_k^{-1/2}m_k+\widetilde Z_k\right)
\mid
\mathcal F_{k-1}
\biggr].
\end{align*}
Moreover,
\[
    \widetilde Z_k
    \sim
    N(0,\Lambda_k),
    \qquad
    \Lambda_k:=V_k^{-1/2}V_{k+1}V_k^{-1/2}.
\]
Since
\[
    V_k=E(\Sigma_k)+V_{k+1},
\]
we have
\[
    I-\Lambda_k
    =
    V_k^{-1/2}E(\Sigma_k)V_k^{-1/2}.
\]

We rewrite the preceding conditional expectation as
\begin{align*}
& E\biggl[
\operatorname{Tr}\left(
\Lambda_k
\nabla^2 f_k\left(V_k^{-1/2}m_k+\widetilde Z_k\right)
\right)
-
\widetilde Z_k^\top
\nabla f_k\left(V_k^{-1/2}m_k+\widetilde Z_k\right)
\mid
\mathcal F_{k-1}
\biggr]
\\
&+
E\biggl[
\operatorname{Tr}\left(
(I-\Lambda_k)
\nabla^2 f_k\left(V_k^{-1/2}m_k+\widetilde Z_k\right)
\right)
\mid
\mathcal F_{k-1}
\biggr]
\\
&-
E\left[
(V_k^{-1/2}m_k)^\top\nabla f_k(\widetilde Z_k)
\mid
\mathcal F_{k-1}
\right]
\\
&-
E\biggl[
(V_k^{-1/2}m_k)^\top
\left(
\nabla f_k\left(V_k^{-1/2}m_k+\widetilde Z_k\right)
-
\nabla f_k(\widetilde Z_k)
\right)
\mid
\mathcal F_{k-1}
\biggr].
\end{align*}
The first term is zero by Gaussian integration by parts: conditional on
$\mathcal F_k$, the quantity $V_k^{-1/2}m_k$ is fixed and
$\widetilde Z_k\sim N(0,\Lambda_k)$ is independent of $\mathcal F_k$. The third term is zero because $f_k$ is $\mathcal F_{k-1}$-measurable as a
random function, $\widetilde Z_k$ is independent of $\mathcal F_k$, and
$E[m_k\mid\mathcal F_{k-1}]=0$.

Using the fundamental theorem of calculus, we obtain
\begin{align*}
& E\left[
h(S_k+T_{k+1})-h(S_{k-1}+T_k)
\mid
\mathcal F_{k-1}
\right]
\\
&=
E\biggl[
\operatorname{Tr}\left(
V_k^{-1/2}E(\Sigma_k)V_k^{-1/2}
\nabla^2 f_k\left(V_k^{-1/2}m_k+\widetilde Z_k\right)
\right)
\mid
\mathcal F_{k-1}
\biggr]
\\
&\quad
-
E\biggl[
(V_k^{-1/2}m_k)^\top
\left(
\nabla^2 f_k(\Theta V_k^{-1/2}m_k+\widetilde Z_k)
-
\nabla^2 f_k(\widetilde Z_k)
\right)
(V_k^{-1/2}m_k)
\mid
\mathcal F_{k-1}
\biggr]
\\
&\quad
-
E\biggl[
(V_k^{-1/2}m_k)^\top
\nabla^2 f_k(\widetilde Z_k)
(V_k^{-1/2}m_k)
\mid
\mathcal F_{k-1}
\biggr],
\end{align*}
where \(\Theta\sim {\rm Uniform}[0,1]\) is independent of all other random
variables.

Define
\[
    A_{k,h}(\mathcal F_{k-1})
    :=
    E\left[
    \nabla^2 f_k(\widetilde Z_k)
    \mid
    \mathcal F_{k-1}
    \right].
\]
Using
\[
    \Sigma_k=E[m_km_k^\top\mid\mathcal F_{k-1}],
\]
the last term becomes
\[
    -
    \operatorname{Tr}\left(
    V_k^{-1/2}\Sigma_k V_k^{-1/2}
    A_{k,h}(\mathcal F_{k-1})
    \right).
\]
By statement (1) of Lemma~\ref{lem: regularity}, the second term above is
bounded in absolute value by
\[
\tilde C_1(d,\beta)\|V_k^{1/2}\|_{op}
E\left[
\|V_k^{-1/2}m_k\|^{2+\beta}
\mid
\mathcal F_{k-1}
\right].
\]
Therefore,
\begin{align*}
& E\left[
h(S_k+T_{k+1})-h(S_{k-1}+T_k)
\mid
\mathcal F_{k-1}
\right]
\\
&\le
E\biggl[
\operatorname{Tr}\left(
V_k^{-1/2}E(\Sigma_k)V_k^{-1/2}
\nabla^2 f_k\left(V_k^{-1/2}m_k+\widetilde Z_k\right)
\right)
\mid
\mathcal F_{k-1}
\biggr]
\\
&\quad
-
\operatorname{Tr}\left(
V_k^{-1/2}\Sigma_k V_k^{-1/2}
A_{k,h}(\mathcal F_{k-1})
\right)
\\
&\quad
+
\tilde C_1(d,\beta)\|V_k^{1/2}\|_{op}
E\left[
\|V_k^{-1/2}m_k\|^{2+\beta}
\mid
\mathcal F_{k-1}
\right].
\end{align*}

Next, add and subtract
\[
\operatorname{Tr}\left(
V_k^{-1/2}E(\Sigma_k)V_k^{-1/2}
A_{k,h}(\mathcal F_{k-1})
\right).
\]
This gives
\begin{align*}
& E\left[
h(S_k+T_{k+1})-h(S_{k-1}+T_k)
\mid
\mathcal F_{k-1}
\right]
\\
&\le
\operatorname{Tr}\left(
A_{k,h}(\mathcal F_{k-1})
V_k^{-1/2}
\left(E(\Sigma_k)-\Sigma_k\right)
V_k^{-1/2}
\right)
\\
&\quad
+
E\biggl[
\operatorname{Tr}\left(
D_k
\left[
\nabla^2 f_k\left(V_k^{-1/2}m_k+\widetilde Z_k\right)
-
\nabla^2 f_k(\widetilde Z_k)
\right]
\right)
\mid
\mathcal F_{k-1}
\biggr]
\\
&\quad
+
\tilde C_1(d,\beta)\|V_k^{1/2}\|_{op}
E\left[
\|V_k^{-1/2}m_k\|^{2+\beta}
\mid
\mathcal F_{k-1}
\right],
\end{align*}
where
\[
    D_k:=V_k^{-1/2}E(\Sigma_k)V_k^{-1/2}.
\]

It remains to bound the weighted trace term. Since
\[
    V_k=E(\Sigma_k)+V_{k+1},
\]
we have
\[
    D_k=I-V_k^{-1/2}V_{k+1}V_k^{-1/2}.
\]
Hence
\[
    0\preceq D_k\preceq I.
\]
Therefore, for any \(x,u\in\mathbb R^d\),
\begin{align*}
&\left|
\operatorname{Tr}\left(
D_k
\left[
\nabla^2 f_k(x+u)-\nabla^2 f_k(x)
\right]
\right)
\right|
\\
&\le
\operatorname{Tr}(D_k)
\left\|
\nabla^2 f_k(x+u)-\nabla^2 f_k(x)
\right\|_{op}
\\
&\le
   \operatorname{Tr}(D_k)\,\tilde C_1(d,\beta)\|V_k^{1/2}\|_{op}\|u\|^\beta.
\end{align*}
Applying this bound with
\[
    x=\widetilde Z_k,
    \qquad
    u=V_k^{-1/2}m_k,
\]
and taking conditional expectation gives
\begin{align*}
&\left|
E\biggl[
\operatorname{Tr}\left(
D_k
\left[
\nabla^2 f_k\left(V_k^{-1/2}m_k+\widetilde Z_k\right)
-
\nabla^2 f_k(\widetilde Z_k)
\right]
\right)
\mid
\mathcal F_{k-1}
\biggr]
\right|
\\
&\le
\tilde C_1(d,\beta)\operatorname{Tr}(D_k)\|V_k^{1/2}\|_{op}
E\left[
\|V_k^{-1/2}m_k\|^\beta
\mid
\mathcal F_{k-1}
\right].
\end{align*}

Thus,
\begin{align*}
& E\left[
h(S_k+T_{k+1})-h(S_{k-1}+T_k)
\mid
\mathcal F_{k-1}
\right]
\\
&\le
\operatorname{Tr}\left(
A_{k,h}(\mathcal F_{k-1})
V_k^{-1/2}
\left(E(\Sigma_k)-\Sigma_k\right)
V_k^{-1/2}
\right)
\\
&\quad
+
\tilde C_1(d,\beta)\operatorname{Tr}(D_k)\|V_k^{1/2}\|_{op}
E\left[
\|V_k^{-1/2}m_k\|^\beta
\mid
\mathcal F_{k-1}
\right]
\\
&\quad
+
\tilde C_1(d,\beta)\|V_k^{1/2}\|_{op}
E\left[
\|V_k^{-1/2}m_k\|^{2+\beta}
\mid
\mathcal F_{k-1}
\right].
\end{align*}

Taking expectations, summing over \(k\), substituting into the Lindeberg
decomposition, and using the rescaling identity for the Wasserstein distance,
we obtain, for the fixed \(h\in Lip_1\),
\begin{align*}
& E\left[
h(W_n)-h(n^{-1/2}V_1^{1/2}Z)
\right]
\\
&\le
\frac{1}{\sqrt n}\sum_{k=1}^n
\Biggl[
\tilde C_1(d,\beta)\|V_k^{1/2}\|_{op}
E\left(\|V_k^{-1/2}m_k\|^{2+\beta}\right)
\\
&\qquad\qquad
+
\tilde C_1(d,\beta)\operatorname{Tr}(D_k)\|V_k^{1/2}\|_{op}
E\left(\|V_k^{-1/2}m_k\|^\beta\right)
\Biggr]
\\
&\quad
+
\frac{1}{\sqrt n}
\left|\sum_{k=1}^n
E\left[
\operatorname{Tr}\left(
A_{k,h}(\mathcal F_{k-1})
V_k^{-1/2}
\left(E(\Sigma_k)-\Sigma_k\right)
V_k^{-1/2}
\right)
\right]
\right|.
\end{align*}
Taking the supremum over \(h\in Lip_1\) proves the claimed bound.

Finally, statement (3) of Lemma~\ref{lem: regularity} gives
\[
    \|A_{k,h}(\mathcal F_{k-1})\|_{op}
    \le
    C\|V_k^{1/2}\|_{op}.
\]
This completes the proof.
\end{proof}

We note that, if
\[
\Sigma_k=\Sigma_\infty
\qquad\text{a.s. for all }k,
\]
then $\Delta_n=0$, and the martingale bound reduces to the two Taylor
remainder terms.

\section{Markov Chain Central Limit Theorem}\label{sec: Markov}
Consider a time-homogeneous Markov chain $\{X_k\}_{k\geq 0}$ on a state space $\mathcal{S}$ with initial distribution $\eta$. Let $r(X_k)\in\Re^d$ denote a vector reward and let $\bar r=\pi r$ whenever the chain has invariant distribution $\pi$. We are interested in rates of convergence to normality of the centered additive functional
\[
    \frac{1}{\sqrt n}\sum_{k=1}^n\bigl(r(X_{k-1})-\bar r\bigr).
\]
As is common in the literature, we study this scaled sum by relating the reward process to a martingale via Poisson's equation \citep{douc2018markov,makowski2002poisson,glynn1996liapounov,glynn2023solution}.

\subsection{Finite State Space Markov Chains}
\label{subsec: finite_state_MC_corrected}

We now show that the martingale theorem can be used to obtain a
non-asymptotic CLT for finite-state Markov chains.  In Theorem~\ref{thm: martingale CLT}, one must control
the signed covariance-interaction term.  This can
be done by applying a second Poisson equation, now to the conditional
covariance function.

Throughout this subsection, assume that the Markov chain is stationary,
finite-state, and irreducible.  Periodicity is allowed.  Let
\[
    g(x):=r(x)-\bar r,
    \qquad
    \bar r:=\pi r.
\]
Let \(V\) solve the Poisson equation
\[
    V-PV=g.
\]
Then
\[
    g(X_{k-1})
    =
    V(X_{k-1})-V(X_k)+m_k,
\]
where
\[
    m_k:=V(X_k)-E[V(X_k)\mid X_{k-1}]
\]
is a martingale difference sequence.  Define
\[
    \Sigma_k:=E[m_km_k^\top\mid \mathcal F_{k-1}]
    =
    v(X_{k-1}),
\]
where
\[
    v(x):=E[m_1m_1^\top\mid X_0=x].
\]
Since the chain is stationary,
\[
    E(\Sigma_k)=\pi(v)=:\Sigma_\infty.
\]
Assume throughout this subsection that \(\Sigma_\infty\) is positive definite.

We first state the technical estimate needed to control the covariance
interaction term.

\begin{lemma}
\label{lem: H_variation}
Let
\[
    V_k:=\sum_{j=k}^n E(\Sigma_j)=(n-k+1)\Sigma_\infty,
    \qquad k=1,\ldots,n,
\]
and
\[
    H_k:=V_k^{-1/2}A_{k,h}(\mathcal F_{k-1})V_k^{-1/2},
\]
where
\[
    A_{k,h}(\mathcal F_{k-1})
    :=E[\nabla^2f_{k,h}(\widetilde Z_k)\mid\mathcal F_{k-1}]
\]
is the matrix appearing in Theorem~\ref{thm: martingale CLT}. Then, for every
$\beta\in(0,1)$,
\[
    \sup_{h\in Lip_1}\sum_{k=2}^nE\|H_k-H_{k-1}\|_{op}
    \le \frac{C}{(1-\beta)^2}n^{(1-\beta)/2},
\]
where $C$ depends on the chain, $r$, $d$, and $\Sigma_\infty$, but not on
$n$ or $\beta$.
\end{lemma}

\begin{proof}
Fix $h\in Lip_1$. Let $c_t:=\sqrt{1-e^{-2t}}$, and let
\[
    P_Ah(x):=E[h(x+A^{1/2}Z)],\qquad Z\sim N(0,I),
\]
denote Gaussian convolution with covariance $A$. The Mehler representation
of the standard normal Stein solution, followed by differentiation under the
integral, gives
\begin{equation}\label{eq: Hk_Mehler_representation}
    H_k=-\int_0^\infty e^{-2t}
    E_{\xi\sim N(0,V_{k+1})}
    [\nabla^2P_{c_t^2V_k}h(S_{k-1}+e^{-t}\xi)]\,dt.
\end{equation}
The differentiation is justified by the standard Gaussian-smoothing bounds
\[
    \|\nabla^jP_Ah\|_{op}
    \le C\|A^{-1/2}\|_{op}^{j-1},\qquad j=2,3,4.
\]
Here $\|\cdot\|_{op}$ denotes the usual matrix operator norm for matrices and Hessians, and the multilinear operator norm for high-order derivatives.

The second- and third-derivative bounds imply, respectively,
\[
\begin{split}
\|\nabla^2P_Ah(x)-\nabla^2P_Ah(y)\|_{op}
&\le 2C\|A^{-1/2}\|_{op},\\
\|\nabla^2P_Ah(x)-\nabla^2P_Ah(y)\|_{op}
&\le C\|A^{-1/2}\|_{op}^2\|x-y\|.
\end{split}
\]
Since $\min\{1,u\}\le u^\beta$ for $u\ge0$, interpolation gives
\begin{equation}\label{eq: uniform_fractional_smoothing}
\|\nabla^2P_Ah(x)-\nabla^2P_Ah(y)\|_{op}
\le C\|A^{-1/2}\|_{op}^{1+\beta}\|x-y\|^\beta,
\end{equation}
with $C$ uniform over $\beta\in(0,1)$.

Because the state space is finite, $m_k$ is uniformly bounded. Moreover,
\[
    c(n-k+1)I\preceq V_k\preceq C(n-k+1)I.
\]
Using \eqref{eq: Hk_Mehler_representation}, compare $H_k$ and $H_{k-1}$ by
changing separately the spatial shift, the covariance inside the Gaussian
convolution, and the covariance of the outer Gaussian average. Denote the
corresponding contributions by $\Delta_{{\rm shift},k}$,
$\Delta_{{\rm int},k}$, and $\Delta_{{\rm ext},k}$.

Since $S_{k-1}-S_{k-2}=m_{k-1}$, \eqref{eq: uniform_fractional_smoothing}
gives
\begin{align*}
E\|\Delta_{{\rm shift},k}\|_{op}
&\le C E\|m_{k-1}\|^\beta
\int_0^\infty e^{-2t}\|(c_t^2V_k)^{-1/2}\|_{op}^{1+\beta}dt\\
&\le \frac{C}{1-\beta}(n-k+1)^{-(1+\beta)/2},
\end{align*}
where
\[
    \int_0^\infty e^{-2t}c_t^{-(1+\beta)}dt=\frac1{1-\beta}.
\]

For the internal covariance change, stationarity gives
$V_{k-1}-V_k=\Sigma_\infty$. Along
$A_s=c_t^2(V_k+s\Sigma_\infty)$, the Gaussian-semigroup identity and the
fourth-derivative bound yield
\begin{align*}
E\|\Delta_{{\rm int},k}\|_{op}
&\le C\int_0^\infty e^{-2t}c_t^2
\|(c_t^2V_k)^{-1/2}\|_{op}^3\|\Sigma_\infty\|_{op}\,dt\\
&\le C(n-k+1)^{-3/2}\int_0^\infty e^{-2t}c_t^{-1}dt\\
&\le C(n-k+1)^{-3/2}.
\end{align*}

For the external covariance change, couple
\[
\xi_{k-1}=\xi_k+\eta_k,\qquad
\xi_k\sim N(0,V_{k+1}),\quad
\eta_k\sim N(0,\Sigma_\infty),
\]
with independent summands. Applying
\eqref{eq: uniform_fractional_smoothing} to the additional shift
$e^{-t}\eta_k$ gives
\[
E\|\Delta_{{\rm ext},k}\|_{op}
\le \frac{C}{1-\beta}(n-k+1)^{-(1+\beta)/2}.
\]
Consequently,
\[
E\|H_k-H_{k-1}\|_{op}
\le \frac{C}{1-\beta}(n-k+1)^{-(1+\beta)/2}
+C(n-k+1)^{-3/2}.
\]
Finally,
\[
\sum_{q=1}^{n-1}q^{-(1+\beta)/2}
\le \frac{C}{1-\beta}n^{(1-\beta)/2},
\]
whereas $\sum_{q\ge1}q^{-3/2}<\infty$. Summing over $k$ and taking the
supremum over $h\in Lip_1$ proves the result.
\end{proof}

We can now prove the finite-state Markov-chain CLT.

\begin{theorem}\label{thm: MC CLT}
Let
\[
U_n=\frac1{\sqrt n}\sum_{k=1}^n(r(X_{k-1})-\bar r).
\]
Assume that the Markov chain is stationary, finite-state, and irreducible,
and that $\Sigma_\infty\succ0$. Then, for every $\beta\in(0,1)$,
\[
d_{\mathcal W}(U_n,\Sigma_\infty^{1/2}Z)
\le \frac{C}{(1-\beta)^2}n^{-\beta/2},
\]
where $C$ does not depend on $n$ or $\beta$. Consequently,
\[
d_{\mathcal W}(U_n,\Sigma_\infty^{1/2}Z)
=O\left(\frac{(\log n)^2}{\sqrt n}\right).
\]
\end{theorem}

\begin{proof}
The reward Poisson equation gives
\[
U_n=W_n+\frac{V(X_0)-V(X_n)}{\sqrt n},
\qquad
W_n:=\frac1{\sqrt n}\sum_{k=1}^n m_k.
\]
Since $V$ is bounded,
\[
d_{\mathcal W}(U_n,\Sigma_\infty^{1/2}Z)
\leq
d_{\mathcal W}(W_n,\Sigma_\infty^{1/2}Z)+O(n^{-1/2}).
\]
In stationarity,
\[
E(\Sigma_k)=\Sigma_\infty,
\qquad
V_k=(n-k+1)\Sigma_\infty,
\]
so the Gaussian target in Theorem~\ref{thm: martingale CLT} is exactly
$\Sigma_\infty^{1/2}Z$.

The martingale increments are uniformly bounded, and their moments can be
bounded uniformly over $\beta\in(0,1)$. Moreover,
\[
\tilde C_1(d,\beta)\leq \frac{C}{1-\beta}
\quad\text{and}\quad
\sum_{q=1}^n q^{-(1+\beta)/2}
\leq \frac{C}{1-\beta}n^{(1-\beta)/2}.
\]
It follows that the first Taylor remainder in
Theorem~\ref{thm: martingale CLT} is at most
\[
C(1-\beta)^{-2}n^{-\beta/2}.
\]
For the second Taylor remainder,
\[
D_k=\frac1{n-k+1}I,
\qquad
\operatorname{Tr}(D_k)=\frac{d}{n-k+1},
\]
and the same bound follows.

It remains to control $\Delta_n$. Let
\[
B(x):=\Sigma_\infty-v(x).
\]
Since $\pi(B)=0$, there is a bounded matrix-valued solution $\Phi$ of
\[
\Phi-P\Phi=B.
\]
We may normalize this solution so that $\pi(\Phi)=0$. Define
\[
\Gamma_k
:=
\Phi(X_k)-E[\Phi(X_k)\mid\mathcal F_{k-1}].
\]
Then
\[
E[\Gamma_k\mid\mathcal F_{k-1}]=0,
\]
and the Markov property and the Poisson equation give
\begin{align*}
B(X_{k-1})
&=
\Phi(X_{k-1})-P\Phi(X_{k-1})
\\
&=
\Phi(X_{k-1})-\Phi(X_k)+\Gamma_k.
\end{align*}
Because the chain is stationary,
\[
E(\Sigma_k)=\Sigma_\infty,
\qquad
\Sigma_k=v(X_{k-1}).
\]
Hence
\begin{equation}
\label{eq: stationary_covariance_poisson_decomposition}
E(\Sigma_k)-\Sigma_k
=
\Phi(X_{k-1})-\Phi(X_k)+\Gamma_k.
\end{equation}

For fixed $h\in Lip_1$, let
\[
H_k
:=
V_k^{-1/2}
A_{k,h}(\mathcal F_{k-1})
V_k^{-1/2}.
\]
Since $H_k$ is $\mathcal F_{k-1}$-measurable,
\[
E\operatorname{Tr}(H_k\Gamma_k)=0.
\]
Using
\eqref{eq: stationary_covariance_poisson_decomposition}
and summation by parts,
\begin{align*}
&\sum_{k=1}^n
E\operatorname{Tr}
\left(
H_k(E(\Sigma_k)-\Sigma_k)
\right)
\\
&=
E\operatorname{Tr}(H_1\Phi(X_0))
-
E\operatorname{Tr}(H_n\Phi(X_n))
\\
&\quad+
\sum_{k=2}^n
E\operatorname{Tr}
\left(
(H_k-H_{k-1})\Phi(X_{k-1})
\right).
\end{align*}
The bound on $A_{k,h}$ in
Theorem~\ref{thm: martingale CLT} and the identity
$V_k=(n-k+1)\Sigma_\infty$ imply
\[
\|H_k\|_{op}
\leq
C(n-k+1)^{-1/2}.
\]
Since $\Phi$ is bounded, the two boundary terms are $O(1)$,
uniformly over $h\in Lip_1$. Lemma~\ref{lem: H_variation} gives
\[
\sup_{h\in Lip_1}
\sum_{k=2}^n
E\|H_k-H_{k-1}\|_{op}
\leq
\frac{C}{(1-\beta)^2}n^{(1-\beta)/2}.
\]
Consequently,
\begin{align*}
\Delta_n
&\leq
\frac{C}{\sqrt n}
\left(
1+
\frac{1}{(1-\beta)^2}n^{(1-\beta)/2}
\right)
\\
&\leq
\frac{C}{(1-\beta)^2}n^{-\beta/2}.
\end{align*}

Combining the reward boundary estimate, the two Taylor-remainder
estimates, and the bound on $\Delta_n$ proves the first assertion.
Choosing
\[
\beta=1-\frac{2}{\log n}
\]
for all sufficiently large $n$ proves the second.
\end{proof}

The stationarity assumption can be removed under aperiodicity. We record the
extension separately because it is used in the TD-learning application.

\begin{corollary}\label{cor: MC CLT nonstationary}
Assume that the Markov chain is finite-state, irreducible, and aperiodic, with
an arbitrary initial distribution. If $\Sigma_\infty\succ0$, then, for every
$\beta\in(0,1)$,
\[
d_{\mathcal W}(U_n,\Sigma_\infty^{1/2}Z)
\le \frac{C}{(1-\beta)^2}n^{-\beta/2},
\]
with $C$ independent of $n$ and $\beta$, and consequently
\[
d_{\mathcal W}(U_n,\Sigma_\infty^{1/2}Z)
=O\left(\frac{(\log n)^2}{\sqrt n}\right).
\]
\end{corollary}

\begin{proof}
Retain the notation from the stationary theorem. The reward Poisson
equation again gives
\[
U_n=W_n+\frac{V(X_0)-V(X_n)}{\sqrt n},
\qquad
W_n:=\frac1{\sqrt n}\sum_{k=1}^n m_k,
\]
and boundedness of $V$ implies
\[
d_{\mathcal W}(U_n,\Sigma_\infty^{1/2}Z)
\leq
d_{\mathcal W}(W_n,\Sigma_\infty^{1/2}Z)+O(n^{-1/2}).
\]

Let
\[
V_k:=\sum_{j=k}^nE(\Sigma_j),
\qquad
D_k:=V_k^{-1/2}E(\Sigma_k)V_k^{-1/2}.
\]
Geometric ergodicity gives
\[
\|E(\Sigma_k)-\Sigma_\infty\|_{op}
\leq
C\rho^{k-1}
\]
for some $\rho\in(0,1)$. Therefore,
\[
\frac1nV_1
=
\Sigma_\infty+O(n^{-1}).
\]
Since $\Sigma_\infty\succ0$, the matrix square-root map is locally
Lipschitz on the positive-definite cone. Coupling the two Gaussian
vectors with the same standard normal vector therefore yields
\[
d_{\mathcal W}
\left(
n^{-1/2}V_1^{1/2}Z,
\Sigma_\infty^{1/2}Z
\right)
=
O(n^{-1}).
\]

Geometric convergence to the positive-definite matrix
$\Sigma_\infty$ also implies that, for all sufficiently large $n$,
there are constants $0<c<C<\infty$ such that
\begin{equation}
\label{eq: nonstationary_covariance_sum_bounds}
c(n-k+1)I
\preceq
V_k
\preceq
C(n-k+1)I,
\qquad
k=1,\ldots,n.
\end{equation}
Indeed, $E(\Sigma_j)\succeq
\frac12\lambda_{\min}(\Sigma_\infty)I$ for all sufficiently large
$j$, while the matrices $E(\Sigma_j)$ are uniformly bounded. For
indices $k$ preceding this fixed threshold, the positive-definite
tail of the sum still contains a fixed positive fraction of the
$n-k+1$ terms.

It follows from
\eqref{eq: nonstationary_covariance_sum_bounds} that
\[
\|V_k^{1/2}\|_{op}
\leq
C(n-k+1)^{1/2},
\qquad
\|V_k^{-1/2}\|_{op}
\leq
C(n-k+1)^{-1/2},
\]
and
\[
\operatorname{Tr}(D_k)
=
\operatorname{Tr}(V_k^{-1}E(\Sigma_k))
\leq
\frac{C}{n-k+1}.
\]
Since the martingale increments are uniformly bounded, the two
Taylor remainders in Theorem~\ref{thm: martingale CLT} are therefore
bounded by
\[
\frac{C}{(1-\beta)^2}n^{-\beta/2}.
\]

It remains to control the covariance-interaction term. Let
\[
B(x):=\Sigma_\infty-v(x),
\]
and let $\Phi$ be a bounded matrix-valued solution of
\[
\Phi-P\Phi=B.
\]
Define
\[
\widetilde\Phi_k
:=
\Phi(X_k)-E[\Phi(X_k)],
\qquad
\Gamma_k
:=
\Phi(X_k)-E[\Phi(X_k)\mid\mathcal F_{k-1}].
\]
The Poisson equation gives
\[
B(X_{k-1})
=
\Phi(X_{k-1})-\Phi(X_k)+\Gamma_k.
\]
Taking expectations and subtracting the resulting identity gives
\[
B(X_{k-1})-E[B(X_{k-1})]
=
\widetilde\Phi_{k-1}
-
\widetilde\Phi_k
+
\Gamma_k.
\]
On the other hand, since $\Sigma_k=v(X_{k-1})$,
\begin{align*}
E(\Sigma_k)-\Sigma_k
&=
E[v(X_{k-1})]-v(X_{k-1})
\\
&=
B(X_{k-1})-E[B(X_{k-1})].
\end{align*}
Therefore,
\begin{equation}
\label{eq: nonstationary_covariance_poisson_decomposition}
E(\Sigma_k)-\Sigma_k
=
\widetilde\Phi_{k-1}
-
\widetilde\Phi_k
+
\Gamma_k.
\end{equation}
Thus, unlike in the stationary theorem, centering is needed because
$E[\Phi(X_k)]$ may depend on $k$.

For fixed $h\in Lip_1$, define
\[
H_k
:=
V_k^{-1/2}
A_{k,h}(\mathcal F_{k-1})
V_k^{-1/2}.
\]
Because $H_k$ is $\mathcal F_{k-1}$-measurable,
\[
E\operatorname{Tr}(H_k\Gamma_k)=0.
\]
Using
\eqref{eq: nonstationary_covariance_poisson_decomposition}
and summation by parts gives
\begin{align*}
&\sum_{k=1}^n
E\operatorname{Tr}
\left(
H_k(E(\Sigma_k)-\Sigma_k)
\right)
\\
&=
E\operatorname{Tr}(H_1\widetilde\Phi_0)
-
E\operatorname{Tr}(H_n\widetilde\Phi_n)
\\
&\quad+
\sum_{k=2}^n
E\operatorname{Tr}
\left(
(H_k-H_{k-1})\widetilde\Phi_{k-1}
\right).
\end{align*}
The matrices $\widetilde\Phi_k$ are uniformly bounded. Moreover,
\eqref{eq: nonstationary_covariance_sum_bounds}
and the bound on $A_{k,h}$ imply
\[
\|H_k\|_{op}
\leq
C(n-k+1)^{-1/2},
\]
so the two boundary terms are $O(1)$.

The proof of Lemma~\ref{lem: H_variation} uses stationarity only
through the linear covariance-sum bounds and the bounded one-step
covariance increments. Here these ingredients are supplied by
\eqref{eq: nonstationary_covariance_sum_bounds} together with
\[
V_{k-1}-V_k=E(\Sigma_{k-1}),
\qquad
V_k-V_{k+1}=E(\Sigma_k),
\]
whose right-hand sides are uniformly bounded positive-semidefinite
matrices. The same spatial-shift, internal-covariance, and
external-covariance estimates therefore give
\[
\sup_{h\in Lip_1}
\sum_{k=2}^n
E\|H_k-H_{k-1}\|_{op}
\leq
\frac{C}{(1-\beta)^2}n^{(1-\beta)/2}.
\]
Consequently,
\[
\Delta_n
\leq
\frac{C}{(1-\beta)^2}n^{-\beta/2}.
\]

Applying Theorem~\ref{thm: martingale CLT} to $W_n$, comparing its
Gaussian target with $\Sigma_\infty^{1/2}Z$, and then using the
reward boundary term proves
\[
d_{\mathcal W}(U_n,\Sigma_\infty^{1/2}Z)
\leq
\frac{C}{(1-\beta)^2}n^{-\beta/2}.
\]
Choosing
\[
\beta=1-\frac{2}{\log n}
\]
for all sufficiently large $n$ gives the stated
\[
O\left(\frac{(\log n)^2}{\sqrt n}\right)
\]
rate. Enlarging $C$ if necessary covers the finitely many remaining
values of $n$.
\end{proof}

\subsection{Extension to General State-Space Markov Chains}

We now extend the results to general state-space Markov chains. The ideas are similar to those of the finite state-space case, but we need additional assumptions and work to prove the result. To avoid unnecessary repetition, we use much of the same notation as in the previous section, while additional terminology and machinery are borrowed from \citep{meyn2012markov}. Not surprisingly, we need stronger assumptions than the conditions needed in \citep{meyn2012markov} to establish the CLT since we are interested in obtaining the rate of convergence.

\begin{theorem}\label{thm: general}
Let $X$ be a $\psi$-irreducible, aperiodic Markov chain on a state space
$\mathcal S$, with transition kernel $P$ and invariant distribution $\pi$.
Suppose there exists a Lyapunov function $\ve\geq 1$, unbounded off petite
sets, such that $E[\ve(X_0)]<\infty$ and
\begin{equation}\label{eq: exp drift}
    P\ve(x)
    =
    E[\ve(X_{k+1})\mid X_k=x]
    \leq
    \lambda\ve(x)+L
\end{equation}
for some $\lambda\in[0,1)$ and $L<\infty$. Let $r:\mathcal S\to\Re^d$ satisfy $\pi\|r\|<\infty$,
let $\bar r:=\pi r$, and define
\[
    g:=r-\bar r,
    \qquad
    U_n:=\frac1{\sqrt n}\sum_{k=1}^n g(X_{k-1}).
\]
Then the following statements hold.
\begin{enumerate}
    \item Fix $\beta\in(0,1)$. If, for some $C_r<\infty$,
    \[
        \|g(x)\|^{2+\beta}\leq C_r\ve(x),
        \qquad x\in\mathcal S,
    \]
    then the Poisson equation $V_r-PV_r=g$ has a solution with the growth
    bound used below. Define
    \[
        m_k:=V_r(X_k)-E[V_r(X_k)\mid X_{k-1}],
        \qquad
        \Sigma_\infty:=E_\pi[m_1m_1^\top].
    \]
    If $\Sigma_\infty\succ0$, then
    \[
        d_{\mathcal W}(U_n,\Sigma_\infty^{1/2}Z)
        \leq
        \frac{C}{(1-\beta)^2}n^{-\beta/2},
    \]
    where $Z\sim N(0,I)$ and $C$ does not depend on $n$.

    \item If
    \[
        \|g(x)\|^4\leq C_r\ve(x),
        \qquad x\in\mathcal S,
    \]
    define $V_r$, $m_k$, and $\Sigma_\infty$ as above. If
    $\Sigma_\infty\succ0$, then
    \[
        d_{\mathcal W}(U_n,\Sigma_\infty^{1/2}Z)
        =
        \mathcal O\!\left(\frac{(\log n)^2}{\sqrt n}\right).
    \]
\end{enumerate}
\end{theorem}

\begin{proof}
The theorem permits an arbitrary initial distribution satisfying
$E[\ve(X_0)]<\infty$. We therefore follow the nonstationary argument of
Corollary~\ref{cor: MC CLT nonstationary}, replacing boundedness by weighted
moment and weighted Poisson-equation estimates. We first prove the first
statement. Set
\[
    p:=2+\beta,
    \qquad
    a:=\frac1p,
    \qquad
    b:=\frac2p.
\]
For every $s\in(0,1]$, Jensen's inequality and
\eqref{eq: exp drift} give
\[
    P\ve^s(x)
    \leq
    (P\ve(x))^s
    \leq
    \lambda^s\ve(x)^s+L^s.
\]
The geometric-drift and weighted Poisson-equation results in
\cite[Chapters~16--17]{meyn2012markov}, applied componentwise, therefore
imply the following facts. First,
\begin{equation}\label{eq: weighted-geometric-general}
    |P^j f(x)-\pi f|
    \leq
    C_s\rho_s^j\|f\|_{\ve^s}
    \bigl(1+\ve(x)^s\bigr),
    \qquad
    \rho_s\in(0,1),
\end{equation}
for functions with
\[
    \|f\|_{\ve^s}
    :=
    \sup_x\frac{|f(x)|}{1+\ve(x)^s}<\infty.
\]
Second, the centered reward Poisson equation has a solution satisfying
\begin{equation}\label{eq: reward-poisson-growth-general}
    \|V_r(x)\|
    \leq
    C\bigl(1+\ve(x)^a\bigr).
\end{equation}
The same conclusions apply componentwise to vector- and matrix-valued
functions. The drift condition also gives
\begin{equation}\label{eq: uniform-V-moment-general}
    \sup_{j\geq0}E[\ve(X_j)]<\infty.
\end{equation}

The martingale increments have uniformly bounded $p$th moments. Indeed,
using \eqref{eq: reward-poisson-growth-general} and conditional Jensen's
inequality,
\begin{align*}
E\|m_k\|^p
&\leq
C E\left[
    \|V_r(X_k)\|^p
    +
    \|E[V_r(X_k)\mid X_{k-1}]\|^p
\right]\\
&\leq
C E\left[1+\ve(X_k)+P\ve(X_{k-1})\right]
\leq C.
\end{align*}
Moreover, \eqref{eq: reward-poisson-growth-general} and
\eqref{eq: uniform-V-moment-general} imply
\[
    \sup_{j\geq0}E\|V_r(X_j)\|<\infty.
\]
Consequently,
\[
    U_n
    =
    W_n+
    \frac{V_r(X_0)-V_r(X_n)}{\sqrt n},
    \qquad
    W_n:=\frac1{\sqrt n}\sum_{k=1}^n m_k,
\]
and the reward boundary term contributes only $\mathcal O(n^{-1/2})$ in
Wasserstein distance.

Define the conditional covariance function
\[
    v(x):=E[m_1m_1^\top\mid X_0=x],
    \qquad
    \Sigma_k=v(X_{k-1}).
\]
The growth estimate \eqref{eq: reward-poisson-growth-general} and
conditional Jensen's inequality yield
\begin{equation}\label{eq: covariance-growth-general}
    \|v(x)\|_{op}
    \leq
    C\bigl(1+\ve(x)^b\bigr).
\end{equation}
In particular,
\[
    \Sigma_\infty=\pi v
\]
is finite. By \eqref{eq: weighted-geometric-general}, after integrating with
respect to the law of $X_0$ and using $E[\ve(X_0)]<\infty$,
\begin{equation}\label{eq: expected-covariance-geometric-general}
    \|E(\Sigma_k)-\Sigma_\infty\|_{op}
    \leq C\rho^k
\end{equation}
for some $\rho\in(0,1)$.

Let
\[
    V_k:=\sum_{j=k}^nE(\Sigma_j),
    \qquad
    q_k:=n-k+1.
\]
Since $\Sigma_\infty\succ0$, the geometric convergence in
\eqref{eq: expected-covariance-geometric-general} and the uniform upper
bound following from \eqref{eq: covariance-growth-general} imply that, for
all sufficiently large $n$,
\begin{equation}\label{eq: covariance-sum-general}
    c q_k I\preceq V_k\preceq Cq_k I,
    \qquad
    k=1,\ldots,n.
\end{equation}
For the lower bound, choose $J$ such that
\[
    E(\Sigma_j)\succeq
    \tfrac12\lambda_{\min}(\Sigma_\infty)I,
    \qquad j\geq J.
\]
If $k\geq J$, the lower bound follows term by term. If $k<J$, the positive
tail from $J$ to $n$ still contains a fixed positive fraction of the
$q_k$ terms once $n$ is sufficiently large. The upper bound follows from
\eqref{eq: covariance-growth-general} and
\eqref{eq: uniform-V-moment-general}. Enlarging the constants, if necessary,
covers the finitely many remaining values of $n$. Furthermore,
\begin{equation}\label{eq: target-covariance-general}
    \frac1nV_1
    =
    \Sigma_\infty+\mathcal O(n^{-1}).
\end{equation}
Since $\Sigma_\infty\succ0$, the matrix square-root map is locally Lipschitz
on the positive-definite cone. Coupling the two Gaussian vectors with the
same standard normal vector therefore gives
\[
    d_{\mathcal W}
    \left(
        n^{-1/2}V_1^{1/2}Z,
        \Sigma_\infty^{1/2}Z
    \right)
    =
    \mathcal O(n^{-1}).
\]

We now apply Theorem~\ref{thm: martingale CLT} to $W_n$. From
\eqref{eq: covariance-sum-general},
\[
    \|V_k^{1/2}\|_{op}\leq Cq_k^{1/2},
    \qquad
    \|V_k^{-1/2}\|_{op}\leq Cq_k^{-1/2}.
\]
Also,
\[
    \operatorname{Tr}(D_k)
    =
    \operatorname{Tr}\bigl(V_k^{-1}E(\Sigma_k)\bigr)
    \leq
    \frac{C}{q_k},
\]
because $E\|\Sigma_k\|_{op}$ is uniformly bounded. The uniform $p$th
moment bound on $m_k$ gives
\[
    E\|V_k^{-1/2}m_k\|^p\leq Cq_k^{-p/2},
    \qquad
    E\|V_k^{-1/2}m_k\|^\beta\leq Cq_k^{-\beta/2}.
\]
Since $\widetilde C_1(d,\beta)\leq C/(1-\beta)$, both Taylor remainders in
Theorem~\ref{thm: martingale CLT} are bounded by a constant times
\begin{align*}
\frac1{\sqrt n}\frac1{1-\beta}
\sum_{k=1}^n q_k^{-(1+\beta)/2}
&\leq
\frac{C}{(1-\beta)^2}n^{-\beta/2},
\end{align*}
where
\[
    \sum_{q=1}^n q^{-(1+\beta)/2}
    \leq
    \frac{C}{1-\beta}n^{(1-\beta)/2}.
\]

It remains to control the covariance-interaction term. As in Corollary~\ref{cor: MC CLT nonstationary}, define
\[
    B(x):=\Sigma_\infty-v(x).
\]
Then $\pi (B)=0$, and \eqref{eq: covariance-growth-general} implies that the
matrix-valued Poisson equation
\begin{equation}\label{eq: covariance-poisson-general}
    \Psi-P\Psi=B
\end{equation}
has a componentwise solution satisfying
\begin{equation}\label{eq: covariance-poisson-growth-general}
    \|\Psi(x)\|_{op}
    \leq
    C\bigl(1+\ve(x)^b\bigr).
\end{equation}
Define
\[
    \widetilde\Psi_k
    :=
    \Psi(X_k)-E[\Psi(X_k)],
    \qquad
    \Gamma_k
    :=
    \Psi(X_k)-E[\Psi(X_k)\mid\mathcal F_{k-1}].
\]
The time-dependent centering is needed here because the theorem does not
assume that $X_0\sim\pi$. By the Markov property and
\eqref{eq: covariance-poisson-general},
\begin{align*}
B(X_{k-1})
&=
\Psi(X_{k-1})-P\Psi(X_{k-1})\\
&=
\Psi(X_{k-1})-\Psi(X_k)+\Gamma_k.
\end{align*}
Taking expectations and subtracting the resulting identity gives
\[
    B(X_{k-1})-E[B(X_{k-1})]
    =
    \widetilde\Psi_{k-1}-\widetilde\Psi_k+\Gamma_k.
\]
Since $\Sigma_k=v(X_{k-1})$,
\[
    E(\Sigma_k)-\Sigma_k
    =
    B(X_{k-1})-E[B(X_{k-1})].
\]
Therefore,
\begin{equation}\label{eq: covariance-interaction-identity-general}
    E(\Sigma_k)-\Sigma_k
    =
    \widetilde\Psi_{k-1}-\widetilde\Psi_k+\Gamma_k.
\end{equation}

For fixed $h\in Lip_1$, let
\[
    H_k
    :=
    V_k^{-1/2}
    A_{k,h}(\mathcal F_{k-1})
    V_k^{-1/2}.
\]
Because $H_k$ is $\mathcal F_{k-1}$-measurable and
$E[\Gamma_k\mid\mathcal F_{k-1}]=0$,
\[
    E\operatorname{Tr}(H_k\Gamma_k)=0.
\]
Using \eqref{eq: covariance-interaction-identity-general} and summation by
parts gives
\begin{align}
&\sum_{k=1}^n
E\operatorname{Tr}\bigl(H_k(E(\Sigma_k)-\Sigma_k)\bigr)
\notag\\
&\quad=
E\operatorname{Tr}(H_1\widetilde\Psi_0)
-
E\operatorname{Tr}(H_n\widetilde\Psi_n)
\notag\\
&\qquad+
\sum_{k=2}^n
E\operatorname{Tr}
\bigl((H_k-H_{k-1})\widetilde\Psi_{k-1}\bigr).
\label{eq: summation-by-parts-general}
\end{align}
The Stein bound in Theorem~\ref{thm: martingale CLT} and
\eqref{eq: covariance-sum-general} imply
\[
    \|H_k\|_{op}\leq Cq_k^{-1/2}.
\]
Together with \eqref{eq: covariance-poisson-growth-general} and
\eqref{eq: uniform-V-moment-general}, this shows that the two boundary terms
in \eqref{eq: summation-by-parts-general} are $\mathcal O(1)$, uniformly over
$h\in Lip_1$.

For the remaining sum, repeat the three-part comparison in
Lemma~\ref{lem: H_variation}. Write the contributions to
$\|H_k-H_{k-1}\|_{op}$ from the internal smoothing covariance, the outer
Gaussian covariance, and the spatial shift as
$\Delta_{{\rm int},k}$, $\Delta_{{\rm ext},k}$, and
$\Delta_{{\rm shift},k}$, respectively. Since
\[
    V_{k-1}-V_k=E(\Sigma_{k-1}),
    \qquad
    V_k-V_{k+1}=E(\Sigma_k),
\]
and these deterministic covariance increments are uniformly bounded, the
fourth-derivative smoothing estimate gives
\[
E\left[
    \Delta_{{\rm int},k}
    \bigl(1+\ve(X_{k-1})^b\bigr)
\right]
\leq
Cq_k^{-3/2}.
\]
Coupling the outer Gaussian variables by an independent Gaussian increment
with covariance $E(\Sigma_k)$ and using the fractional smoothing estimate
gives
\[
E\left[
    \Delta_{{\rm ext},k}
    \bigl(1+\ve(X_{k-1})^b\bigr)
\right]
\leq
\frac{C}{1-\beta}q_k^{-(1+\beta)/2}.
\]
Finally, since $S_{k-1}-S_{k-2}=m_{k-1}$, the spatial-shift term satisfies
\begin{align*}
&E\left[
    \Delta_{{\rm shift},k}
    \bigl(1+\ve(X_{k-1})^b\bigr)
\right]\\
&\qquad\leq
\frac{C}{1-\beta}q_k^{-(1+\beta)/2}
E\left[
    \|m_{k-1}\|^\beta
    \bigl(1+\ve(X_{k-1})^b\bigr)
\right].
\end{align*}
The final expectation is uniformly bounded. Indeed, since $p-\beta=2$ and
$bp/(p-\beta)=1$, H\"older's inequality gives
\begin{align*}
E\left[\|m_{k-1}\|^\beta\ve(X_{k-1})^b\right]
&\leq
\left(E\|m_{k-1}\|^p\right)^{\beta/p}
\left(E\ve(X_{k-1})^{bp/(p-\beta)}\right)^{(p-\beta)/p}\\
&\leq C.
\end{align*}
Consequently,
\begin{align*}
&\sum_{k=2}^n
E\left[
    \|H_k-H_{k-1}\|_{op}
    \bigl(1+\ve(X_{k-1})^b\bigr)
\right]\\
&\qquad\leq
C\sum_{q=1}^{n-1}q^{-3/2}
+
\frac{C}{1-\beta}
\sum_{q=1}^{n-1}q^{-(1+\beta)/2}\\
&\qquad\leq
\frac{C}{(1-\beta)^2}n^{(1-\beta)/2}.
\end{align*}
By \eqref{eq: covariance-poisson-growth-general}, the same bound controls
the last sum in \eqref{eq: summation-by-parts-general}. After multiplying by
the global $n^{-1/2}$ factor in Theorem~\ref{thm: martingale CLT}, we obtain
\[
    \Delta_n
    \leq
    \frac{C}{(1-\beta)^2}n^{-\beta/2}.
\]
Combining this estimate with the two Taylor remainders, the reward boundary
term, and the Gaussian-target comparison proves the first statement.

For the second statement, the fourth-moment condition gives a reward Poisson
solution with the fixed growth bound
\[
    \|V_r(x)\|\leq C\bigl(1+\ve(x)^{1/4}\bigr),
\]
uniformly bounded fourth moments of $m_k$, and a covariance function $v$ and
covariance Poisson solution $\Psi$ with the fixed growth bound
\[
    \|v(x)\|_{op}+\|\Psi(x)\|_{op}
    \leq
    C\bigl(1+\ve(x)^{1/2}\bigr).
\]
Repeating the preceding argument for any $\beta\in(0,1)$ gives
\[
    d_{\mathcal W}(U_n,\Sigma_\infty^{1/2}Z)
    \leq
    \frac{C}{(1-\beta)^2}n^{-\beta/2},
\]
with $C$ uniform in $\beta$. In particular, the weighted spatial-shift
estimate is uniform because H\"older's inequality gives
\begin{align*}
E\left[\|m_{k-1}\|^\beta\ve(X_{k-1})^{1/2}\right]
&\leq
\left(E\|m_{k-1}\|^4\right)^{\beta/4}
\left(E\ve(X_{k-1})^{2/(4-\beta)}\right)^{(4-\beta)/4},
\end{align*}
and $2/(4-\beta)\leq2/3$. Choosing
\[
    \beta=1-\frac{2}{\log n}
\]
for all sufficiently large $n$ gives
\[
    d_{\mathcal W}(U_n,\Sigma_\infty^{1/2}Z)
    =
    \mathcal O\!\left(\frac{(\log n)^2}{\sqrt n}\right).
\]
Enlarging the implied constant covers the finitely many remaining values of
$n$.
\end{proof}

\section{An Application to TD Learning}\label{sec: TD}

Temporal Difference (TD) Learning is a common method to learn the performance of a policy in reinforcement learning; see \citep{bertsekas1996neuro,sutton2018reinforcement} for details. In \citep{tsitsiklis1996analysis}, it was shown that TD learning can be represented as a jump Markovian linear dynamical system as follows:
\begin{equation}\label{eq: TD}
    \theta_{k+1}=\theta_k-\epsilon_k (A(X_k)\theta_k+b(X_k)),
\end{equation}
where $\theta_k\in \Re^d,$ $\epsilon_k$ is the stepsize, $\{X_k\}_{k\geq 0}$ is a finite state-space Markov chain, $A(X_k)$ is a $d\times d$ matrix, $b(X_k)$ is a $d$-dimensional vector and $\theta_0$ is a fixed vector. We make the following assumptions.
\begin{assumption}\label{assume: TD}
We assume that the TD learning algorithm satisfies the following conditions:
\begin{enumerate}
    \item $\{X_k\}$ is an irreducible, aperiodic, finite state Markov chain. 
    \item Let $\bar{A}:=E_{X_k\sim\pi}(A(X_k)).$ Then, $-\bar{A}$ is a Hurwitz matrix, i.e., all eigenvalues have strictly negative real parts.
\end{enumerate}
\hfill $\diamond$
\end{assumption}
The first item in the assumption is standard in the reinforcement learning literature \citep{sutton2018reinforcement}, while the second item has been established for TD learning in \citep{tsitsiklis1996analysis}.

Define $\theta^*$ by $\bar{A}\theta^*+\bar{b}=0,$ where $\bar{b}=E_{X_k\sim\pi}(b(X_k)).$  Here, we are interested in the Polyak-Ruppert averaged version of the TD-learning algorithm, i.e.,
$$\bar{\theta}_n=\frac{1}{n}\sum_{k=1}^{n} \theta_k.$$ Our main result of this section characterizes the rate of convergence of this algorithm for an appropriately chosen $\{\epsilon_k\}$. As in the previous theorem, we have chosen to write the result using the $\mathcal{O}$ notation rather than making the constants explicit, but these constants can be made explicit by following the details of the proof of the theorem.

Define the centered TD noise
\begin{equation}\label{eq: TD reward}
    r_{\rm TD}(x)
    :=
    (A(x)-\bar A)\theta^*+(b(x)-\bar b).
\end{equation}
Since $\bar A\theta^*+\bar b=0$, we have $\pi r_{\rm TD}=0$. Let
$V_{\rm TD}$ solve
\[
    V_{\rm TD}-PV_{\rm TD}=r_{\rm TD},
\]
and define
\begin{equation}\label{eq: TD asymptotic covariance}
    m_{k+1}
    :=
    V_{\rm TD}(X_{k+1})
    -E[V_{\rm TD}(X_{k+1})\mid X_k],
    \qquad
    \Sigma_\infty
    :=
    E_\pi[m_1m_1^\top].
\end{equation}

\begin{theorem}\label{thm: TD CLT}
Suppose Assumption~\ref{assume: TD} holds, $\Sigma_\infty\succ0$, and
\[
    \epsilon_k=(k+1)^{-\delta},
    \qquad
    \delta\in(1/2,1).
\]
Then, for $Z\sim N(0,I)$,
\[
 d_\mathcal{W}\!\left(
    \sqrt n(\bar\theta_n-\theta^*),
    (\bar A^{-1}\Sigma_\infty\bar A^{-T})^{1/2}Z
 \right)
 =
 \mathcal O\!\left(
    \max\left\{
        \frac{\sqrt{\log n}}{n^{\delta-1/2}},
        \frac1{n^{(1-\delta)/2}}
    \right\}
 \right).
\]
\end{theorem}

\begin{proof}
Let $\Delta_k:=\theta_k-\theta^*$. From \eqref{eq: TD},
\[
\Delta_{k+1}
=(I-\epsilon_k\bar A)\Delta_k
-\epsilon_k(A(X_k)-\bar A)\Delta_k
-\epsilon_k r_{\rm TD}(X_k).
\]
Iterating this recursion, summing over $k=1,\ldots,n$, and interchanging the
order of summation gives
\begin{align*}
\sqrt n(\bar\theta_n-\theta^*)
&=
\frac1{\sqrt n}\sum_{k=1}^n
\prod_{l=0}^{k-1}(I-\epsilon_l\bar A)\Delta_0\\
&\quad-
\frac1{\sqrt n}\sum_{j=0}^{n-1}
\epsilon_j\sum_{k=j+1}^n
\prod_{l=j+1}^{k-1}(I-\epsilon_l\bar A)
(A(X_j)-\bar A)\Delta_j\\
&\quad-
\frac1{\sqrt n}\sum_{j=0}^{n-1}
\epsilon_j\sum_{k=j+1}^n
\prod_{l=j+1}^{k-1}(I-\epsilon_l\bar A)
r_{\rm TD}(X_j),
\end{align*}
where an empty matrix product is the identity. For $0\leq j\leq n$, define
\begin{equation}\label{eq: def Upsilon}
\begin{split}
\Upsilon_j^n
&:=
\epsilon_j\sum_{k=j+1}^{n}
\prod_{l=j+1}^{k-1}(I-\epsilon_l\bar A)
-\bar A^{-1},\\
\Phi_j^n
&:=
\epsilon_j\sum_{k=j}^{n-1}
\prod_{l=j}^{k-1}(I-\epsilon_l\bar A)
-\bar A^{-1},
\qquad 0\leq j\leq n-1.
\end{split}
\end{equation}
Thus $\Upsilon_n^n=-\bar A^{-1}$. Directly from the definitions,
\begin{equation}\label{eq: correct Upsilon Phi identity}
    \Upsilon_j^n
    =
    \frac{\epsilon_j}{\epsilon_{j+1}}\Phi_{j+1}^{n+1}
    +
    \left(\frac{\epsilon_j}{\epsilon_{j+1}}-1\right)\bar A^{-1},
    \qquad 0\leq j\leq n-1.
\end{equation}
Appendix~\ref{sec: Upsilon} proves
\begin{equation}\label{eq: Upsilon-uniform-and-average}
    \sup_{n\geq1}\max_{0\leq j\leq n}\|\Upsilon_j^n\|_{op}<\infty,
    \qquad
    \frac1n\sum_{j=0}^{n-1}\|\Upsilon_j^n\|_{op}^2
    =
    \mathcal O\!\left(n^{-(1-\delta)}\right).
\end{equation}
Using $\epsilon_j\sum\prod=\bar A^{-1}+\Upsilon_j^n$, we obtain the five-term
decomposition
\begin{align}
\sqrt n(\bar\theta_n-\theta^*)
&=
\frac1{\sqrt n}\sum_{k=1}^{n}
\prod_{l=0}^{k-1}(I-\epsilon_l\bar A)\Delta_0
\label{eq1}\\
&\quad-
\frac1{\sqrt n}\sum_{j=0}^{n-1}
\bar A^{-1}(A(X_j)-\bar A)\Delta_j
\label{eq2}\\
&\quad-
\frac1{\sqrt n}\sum_{j=0}^{n-1}
\Upsilon_j^n(A(X_j)-\bar A)\Delta_j
\label{step3}\\
&\quad-
\frac1{\sqrt n}\sum_{j=0}^{n-1}
\Upsilon_j^n r_{\rm TD}(X_j)
\label{step4}\\
&\quad-
\frac1{\sqrt n}\sum_{j=0}^{n-1}
\bar A^{-1}r_{\rm TD}(X_j).
\label{step5}
\end{align}
Denote the five terms on the right-hand side of
\eqref{eq1}--\eqref{step5} by $T_{1,n},\ldots,T_{5,n}$, respectively. We
bound the first four in mean norm and apply the Markov-chain CLT to the fifth.

\paragraph{Step 1: term \eqref{eq1}.}
Because $\epsilon_0=1$ and all factors are polynomials in $\bar A$,
\[
\epsilon_0\sum_{k=1}^n\prod_{l=0}^{k-1}(I-\epsilon_l\bar A)
=
(I-\epsilon_0\bar A)(\Upsilon_0^n+\bar A^{-1}).
\]
The uniform bound in \eqref{eq: Upsilon-uniform-and-average} therefore shows
that the mean norm of \eqref{eq1} is $\mathcal O(n^{-1/2})$.

\paragraph{Step 2: term \eqref{eq2}.}
Let $V_A$ solve the componentwise matrix Poisson equation
\[
    V_A-PV_A=A-\bar A,
\]
and set
\[
    \mu_{j+1}:=V_A(X_{j+1})-E[V_A(X_{j+1})\mid X_j].
\]
The finite state space makes $V_A$ and $\mu_j$ bounded. Ignoring the overall
sign in \eqref{eq2}, Poisson's equation and summation by parts give
\begin{align}
&\frac1{\sqrt n}\sum_{j=0}^{n-1}
\bar A^{-1}(V_A(X_j)\Delta_j-V_A(X_{j+1})\Delta_{j+1})
\label{eq3}\\
&\quad+
\frac1{\sqrt n}\sum_{j=0}^{n-1}
\bar A^{-1}V_A(X_{j+1})(\Delta_{j+1}-\Delta_j)
\label{eq4}\\
&\quad+
\frac1{\sqrt n}\sum_{j=0}^{n-1}
\bar A^{-1}\mu_{j+1}\Delta_j.
\label{eq5}
\end{align}
The first sum telescopes and has mean norm $\mathcal O(n^{-1/2})$. Since the
coefficients $A(x)$ and $b(x)$ are bounded and
$\sup_jE\|\Delta_j\|^2<\infty$,
\[
    E\|\Delta_{j+1}-\Delta_j\|
    \leq C\epsilon_j,
\]
so \eqref{eq4} has mean norm
$\mathcal O(n^{-(\delta-1/2)})$. Finally, martingale orthogonality and
boundedness of $\mu_j$ give
\begin{align*}
E\left\|
\frac1{\sqrt n}\sum_{j=0}^{n-1}
\bar A^{-1}\mu_{j+1}\Delta_j
\right\|^2
&\leq
\frac Cn\sum_{j=0}^{n-1}E\|\Delta_j\|^2\\
&=
\mathcal O\!\left(\frac{\log n}{n^{2\delta-1}}\right),
\end{align*}
where Appendix~\ref{sec: Delta} supplies the last estimate. Hence the term
in \eqref{eq2} satisfies
\begin{equation}\label{eq: TD step2 bound}
    E\|T_{2,n}\|
    =
    \mathcal O\!\left(
        \frac{\sqrt{\log n}}{n^{\delta-1/2}}
    \right).
\end{equation}

\paragraph{Step 3: term \eqref{step3}.}
Using the same matrix Poisson equation,
\begin{align}
&\frac1{\sqrt n}\sum_{j=0}^{n-1}
\bigl(\Upsilon_j^nV_A(X_j)\Delta_j
-\Upsilon_{j+1}^nV_A(X_{j+1})\Delta_{j+1}\bigr)
\label{eq: TD step3 boundary}\\
&\quad+
\frac1{\sqrt n}\sum_{j=0}^{n-1}
\Upsilon_{j+1}^nV_A(X_{j+1})(\Delta_{j+1}-\Delta_j)
\label{eq: TD step3 increment}\\
&\quad+
\frac1{\sqrt n}\sum_{j=0}^{n-1}
(\Upsilon_{j+1}^n-\Upsilon_j^n)V_A(X_{j+1})\Delta_j
\label{eq: TD step3 variation}\\
&\quad+
\frac1{\sqrt n}\sum_{j=0}^{n-1}
\Upsilon_j^n\mu_{j+1}\Delta_j.
\label{eq: TD step3 martingale}
\end{align}
The first term telescopes and is $\mathcal O(n^{-1/2})$ in mean norm. The
second is $\mathcal O(n^{-(\delta-1/2)})$. To handle the third, the definitions
give
\begin{align*}
\Upsilon_{j+1}^n-\Upsilon_j^n
&=
\left(1-\frac{\epsilon_j}{\epsilon_{j+1}}\right)
(\Upsilon_{j+1}^n+\bar A^{-1})\\
&\quad+
\epsilon_j\bar A(\Upsilon_{j+1}^n+\bar A^{-1})
-\epsilon_j I.
\end{align*}
Since
$|1-\epsilon_j/\epsilon_{j+1}|=\mathcal O((j+1)^{-1})
=\mathcal O(\epsilon_j)$ and $\Upsilon_j^n$ is uniformly bounded,
\begin{equation}\label{eq: Upsilon variation TD}
    \|\Upsilon_{j+1}^n-\Upsilon_j^n\|_{op}
    \leq C\epsilon_j.
\end{equation}
Thus \eqref{eq: TD step3 variation} is
$\mathcal O(n^{-(\delta-1/2)})$ in mean norm. Martingale orthogonality treats
\eqref{eq: TD step3 martingale} exactly as \eqref{eq5}, giving the
$\sqrt{\log n}/n^{\delta-1/2}$ rate. Therefore,
\begin{equation}\label{eq: TD step3 bound}
    E\|T_{3,n}\|
    =
    \mathcal O\!\left(
        \frac{\sqrt{\log n}}{n^{\delta-1/2}}
    \right).
\end{equation}

\paragraph{Step 4: term \eqref{step4}.}
The reward Poisson equation gives
\begin{align*}
&\frac1{\sqrt n}\sum_{j=0}^{n-1}\Upsilon_j^n r_{\rm TD}(X_j)\\
&\quad=
\frac1{\sqrt n}\sum_{j=0}^{n-1}
\bigl(\Upsilon_j^nV_{\rm TD}(X_j)
-\Upsilon_{j+1}^nV_{\rm TD}(X_{j+1})\bigr)\\
&\qquad+
\frac1{\sqrt n}\sum_{j=0}^{n-1}
(\Upsilon_{j+1}^n-\Upsilon_j^n)V_{\rm TD}(X_{j+1})\\
&\qquad+
\frac1{\sqrt n}\sum_{j=0}^{n-1}\Upsilon_j^n m_{j+1}.
\end{align*}
The first term is $\mathcal O(n^{-1/2})$ and, by
\eqref{eq: Upsilon variation TD}, the second is
$\mathcal O(n^{-(\delta-1/2)})$. For the martingale term,
\eqref{eq: Upsilon-uniform-and-average} gives
\begin{align*}
E\left\|
\frac1{\sqrt n}\sum_{j=0}^{n-1}\Upsilon_j^n m_{j+1}
\right\|^2
&\leq
\frac Cn\sum_{j=0}^{n-1}\|\Upsilon_j^n\|_{op}^2\\
&=
\mathcal O\!\left(n^{-(1-\delta)}\right).
\end{align*}
Hence
\begin{equation}\label{eq: TD step4 bound}
    E\|T_{4,n}\|
    =
    \mathcal O\!\left(
        n^{-(\delta-1/2)}+n^{-(1-\delta)/2}
    \right).
\end{equation}

\paragraph{Step 5 and conclusion.}
By Corollary~\ref{cor: MC CLT nonstationary} and the fact that a deterministic
linear map is Lipschitz,
\begin{align*}
&d_\mathcal W\!\left(
-
\frac1{\sqrt n}\sum_{j=0}^{n-1}\bar A^{-1}r_{\rm TD}(X_j),
(\bar A^{-1}\Sigma_\infty\bar A^{-T})^{1/2}Z
\right)\\
&\qquad=
\mathcal O\!\left(\frac{(\log n)^2}{\sqrt n}\right).
\end{align*}
For $h\in Lip_1$, successively adding the first four terms in
\eqref{eq1}--\eqref{step4} changes its expectation by at most the sum of their
mean norms. Combining \eqref{eq: TD step2 bound},
\eqref{eq: TD step3 bound}, and \eqref{eq: TD step4 bound} therefore gives
\begin{align*}
&d_\mathcal W\!\left(
\sqrt n(\bar\theta_n-\theta^*),
(\bar A^{-1}\Sigma_\infty\bar A^{-T})^{1/2}Z
\right)\\
&\quad=
\mathcal O\!\left(
\frac1{\sqrt n}
+
\frac{\sqrt{\log n}}{n^{\delta-1/2}}
+
\frac1{n^{(1-\delta)/2}}
+
\frac{(\log n)^2}{\sqrt n}
\right).
\end{align*}
For every fixed $\delta\in(1/2,1)$,
$(\log n)^2/\sqrt n=o(\sqrt{\log n}/n^{\delta-1/2})$. This proves the
stated rate.
\end{proof}

\section{Conclusions}\label{sec: conclusions}

We derived a non-asymptotic central limit theorem for martingales and Markov chains, and applied the results to get rate of convergence in the central limit theorem for TD learning with averaging. The results for martingales and Markov chains may be more broadly applicable to study other reinforcement learning algorithms as well as other algorithms such as nonlinear stochastic approximation and stochastic gradient descent.

There are several possible avenues for further work:
\begin{itemize}
    \item One can try to use the martingale CLT results in this paper to strengthen the results and allow for Markovian noise in \citep{anastasiou2019normal}. However, these extensions would require us to go beyond the linear setting in this paper and consider nonlinear problems. Along similar lines, it would be also interesting to see if our results can be used to study contractive stochastic approximation as in \citep{chen2023lyapunov} or other stochastic approximation schemes \citep{borkar2009stochastic,kushnerg}.
    \item Under some conditions, the rate of convergence in the central limit theorem for i.i.d. random vectors is $\mathcal{O}(1/\sqrt{n})$ \citep{fang2019multivariate}. It would be interesting to determine whether the $\mathcal{O}((\log n)^2/\sqrt{n})$ rate in Theorem~\ref{thm: MC CLT} can be improved to $\mathcal{O}(1/\sqrt{n})$ under additional conditions. The Edgeworth expansion for scalar functions of Markov chains in \citep{kontoyiannis2003spectral} suggests that the logarithmic loss might be removable. In the present proof, one factor of $(1-\beta)^{-1}$ comes from the regularity constant in Stein's equation and a second comes from summing a near-harmonic sequence after the Lindeberg decomposition. A different argument may be required to eliminate these factors.
    \item The approach in \citep{gallouet2018regularity} that we have used is known to be not tight in terms of the dependence on the dimension $d$ for i.i.d. random vectors. Under some conditions, the dependence on $d$ can be considerably strengthened in the i.i.d. case \citep{courtade2019existence,zhai2018high}. Whether that is possible for martingales and Markov chains is an interesting question for further research.
    \item The Wasserstein distance used in this paper is the $1$-Wasserstein distance. It would be interesting to obtain convergence rates in the stronger $p$-Wasserstein distance. Work along these lines under different conditions than in this paper can be found in \citep{bonis2020stein}, and in the previously mentioned references \citep{courtade2019existence,zhai2018high}. It would be interesting to see if such results can be obtained for martingales and Markov chains.
    \item It is interesting to note that, while Theorem~\ref{thm: MC CLT} is critical to obtain convergence to the central limit theorem for TD learning, the rate of convergence is dominated by the parameter $\delta$ used in the algorithm. It would be interesting to see if the convergence rate result is tight or whether it can be improved through other techniques.
    \item Ignoring the polylogarithmic factor in Theorem~\ref{thm: TD CLT}, our result suggests that $2/3$ is a good choice for $\delta$. This is different than the conclusion in \citep{haque2023tight}, where their upper bound implies $\delta=0.75$ is the best choice. However, these bounds ignore the constants hidden in the $\mathcal{O}$ notation, which could depend on $\delta$ (see Appendix~\ref{sec: Upsilon}). Moreover, what we have is only an upper bound, so a more refined analysis is needed to understand a good choice for $\delta.$
\end{itemize}

\paragraph{Acknowledgment:} 
\begin{itemize}
    \item I thank Prof. Siva Theja Maguluri, Georgia Tech, for many stimulating discussions and drawing my attention to several references during the course of this work, Prof. Sean Meyn, University of Florida, for pointing out the bias issue with fixed stepsize TD learning and why it cannot be eliminated by averaging, and other discussions, Prof. Mehrdad Moharrami, University of Iowa, for carefully reading an earlier version of the paper and providing many useful comments, and Weichen Wu, CMU,  for pointing out a missing factor in Theorem 4 in an earlier version of the paper.
    \item I thank Professors Alexey Naumov and Sergey Samsonov for pointing out a gap in the proof of the nonasymptotic martingale CLT in an earlier version of the paper.
    \item Research supported by AFOSR Grant FA9550-24-1-0002, ONR Grant N00014-19-1-2566, and NSF Grants CNS 23-12714, CNS 21-06801, CCF 19-34986, and CCF 22-07547.
\end{itemize}

\bibliographystyle{informs2014}
\bibliography{references}

\appendix

\section{Properties of $E(||\Delta_k||^2)$}\label{sec: Delta}

In this section and the next, constants are local to the appendix; letters used for constants elsewhere in the paper may therefore be reused.
Since the Markov chain is finite-state, irreducible, and aperiodic, its mixing time satisfies $\tau_\epsilon=O(\log(1/\epsilon))$; hence, for $\epsilon_k=(k+1)^{-\delta}$, we have $k-\tau_{\epsilon_k}\geq 0$, $\epsilon_{k-\tau_{\epsilon_k}}/\epsilon_k=O(1)$, and $\epsilon_k\tau_{\epsilon_k}\to0$, so the diminishing-stepsize conditions of \cite[Theorem~12]{srikant2019finite} hold for all sufficiently large $k$. Thus, there exists a positive integer $k^*$ such that for any $k\geq k^*,$ we have
\begin{align*}
E\left[\|{\Delta}_{k}\|^2\right]\leq \kappa_1\left(\prod_{j=k^*}^{k-1} a_j\right) + \kappa_2\sum_{j=k^*}^{k-1}b_j\left(\prod_{l=j+1}^{k-1}a_l\right),
\end{align*}
where $a_j=1-{\kappa_3\epsilon_j},$ $b_j= \kappa_4\epsilon_j^2\log(1/\epsilon_j),$ and $\kappa_1,\kappa_2,\kappa_3,\kappa_4$ are positive constants. Since $a_j\leq \exp(-\kappa_3\epsilon_j),$ the first term is bounded by $C\exp(-c k^{1-\delta})$ and hence goes to zero at a stretched-exponential rate. The second term can be upper bounded as follows:
\begin{align*}
&\kappa_2\kappa_4\sum_{j=0}^{k-1}\epsilon_j^2\log\left(\frac{1}{\epsilon_j}\right)\exp\left(-\kappa_3\sum_{l=j+1}^{k-1}\epsilon_l\right)\\
&\leq \kappa_2\kappa_4\delta \log k\sum_{j=0}^{k-1}\epsilon_j^2\exp\left(-\kappa_3\sum_{l=j+1}^{k-1}\epsilon_l\right)\\
&\leq \kappa_2\kappa_4\delta \log k\left(\sum_{j=0}^{\lceil k/2\rceil}\exp\left(-\kappa_3\sum_{l=j+1}^{k-1}\epsilon_l\right)+\sum_{j=\lceil k/2\rceil+1}^{k-1}\epsilon_j^2\right)\\
&=\kappa_2\kappa_4\delta \log k\left(\sum_{j=0}^{\lceil k/2\rceil}\exp\left(-\kappa_3\sum_{l=j+1}^{k-1}(\frac{1}{l+1})^\delta\right)+\sum_{j=\lceil k/2\rceil+1}^{k-1}(\frac{1}{j+1})^{2\delta}\right)=\mathcal{O}\left(\frac{\log k}{k^{2\delta-1}}\right).
\end{align*}

\section{Properties of $\Upsilon_j^n$ and $\Phi_j^n$}\label{sec: Upsilon}

Recall the definition of $\Upsilon_j^n$ and $\Phi_j^n$ in
(\ref{eq: def Upsilon}). In this section of the appendix, we examine
certain expressions which are shown to be either bounded or converging
to zero in \citep{polyak1992acceleration}, and derive explicit bounds
for them as functions of $\delta$.

From \cite[page 846]{polyak1992acceleration}, we know that $\Phi_j^t$
can be written as
\[
\Phi_j^t=C_j^t+\bar A^{-1}D_j^t,
\]
where
\[
\|C_j^t\|_{op}
\leq
\frac{\epsilon_j-\epsilon_{j+1}}{\epsilon_j}
\sum_{i=j}^t m_j^i e^{-\lambda m_j^i},
\qquad
D_j^t
=
\prod_{l=j}^{t-1}(I-\epsilon_l\bar A),
\]
and
\[
\|D_j^t\|_{op}
\leq
K e^{-\alpha m_j^{t-1}},
\qquad
m_j^i:=\sum_{k=j}^i \epsilon_k,
\]
for $t\geq j+1$ and some constants $K,\alpha,\lambda>0$.

We first bound the contribution of $D_j^t$. Since
$\epsilon_k=(k+1)^{-\delta}$ is nonincreasing, for $0\leq j\leq t-1$,
\[
m_j^{t-1}
=
\sum_{k=j}^{t-1}\epsilon_k
\geq
(t-j)t^{-\delta}.
\]
Hence
\[
\|D_j^t\|_{op}^2
\leq
K^2\exp\!\left(-2\alpha (t-j)t^{-\delta}\right).
\]
Therefore,
\begin{align*}
\sum_{j=0}^{t-1}\|D_j^t\|_{op}^2
&\leq
K^2\sum_{r=1}^{t}\exp\!\left(-2\alpha r t^{-\delta}\right)\\
&\leq
\frac{K^2}{1-\exp(-2\alpha t^{-\delta})}
=
\mathcal O(t^\delta),
\end{align*}
and consequently
\begin{equation}\label{eq: D-average-square}
\frac1t\sum_{j=0}^{t-1}\|D_j^t\|_{op}^2
=
\mathcal O\!\left(\frac1{t^{1-\delta}}\right).
\end{equation}
This rate is natural: for indices $j$ within distance of order
$t^\delta$ from $t$, the product $D_j^t$ contains only a small number of
factors and need not yet be small.

We next bound $C_j^t$. Using
\[
\frac{\epsilon_j-\epsilon_{j+1}}{\epsilon_j}
=
1-\left(\frac{j+1}{j+2}\right)^\delta
\leq
\frac{1}{j+2},
\]
we have
\[
\|C_j^t\|_{op}
\leq
\frac{1}{j+2}\sum_{i=j}^t m_j^i e^{-\lambda m_j^i}.
\]
Fix $\mu>\delta/\lambda$. Since
\[
\epsilon_k=(k+1)^{-\delta}\geq \frac{\mu}{k+1}
\]
for all sufficiently large $k$, there exists an integer
$K_{\mu,\delta}$ such that, for all $i\geq j\geq K_{\mu,\delta}$,
\[
m_j^i
=
\sum_{k=j}^i \epsilon_k
\geq
\mu\sum_{k=j}^i \frac1{k+1}
\geq
\mu\log\frac{i+1}{j+1}.
\]
Thus,
\[
\left(\frac{i+1}{j+1}\right)^\delta
\leq
\exp\!\left(\frac{\delta}{\mu}m_j^i\right).
\]
Using $m_j^i-m_j^{i-1}=\epsilon_i=(i+1)^{-\delta}$, we obtain, for
$j\geq K_{\mu,\delta}$,
\begin{align*}
\|C_j^t\|_{op}
&\leq
\frac{1}{j+2}
\sum_{i=j}^t
m_j^i e^{-\lambda m_j^i}
\frac{m_j^i-m_j^{i-1}}{\epsilon_i}\\
&=
\frac{ (j+1)^\delta}{j+2}
\sum_{i=j}^t
m_j^i e^{-\lambda m_j^i}
\left(\frac{i+1}{j+1}\right)^\delta
(m_j^i-m_j^{i-1})\\
&\leq
\frac{ (j+1)^\delta}{j+2}
\sum_{i=j}^t
m_j^i
\exp\!\left(-\left(\lambda-\frac{\delta}{\mu}\right)m_j^i\right)
(m_j^i-m_j^{i-1})\\
&\leq
\frac{C}{(j+1)^{1-\delta}},
\end{align*}
where $C<\infty$ depends on $\delta,\lambda,\mu$, but not on $j$ or
$t$. Enlarging the constant to absorb the finitely many indices
$j<K_{\mu,\delta}$, we obtain
\begin{equation}\label{eq: C-bound}
\|C_j^t\|_{op}
\leq
\frac{C}{(j+1)^{1-\delta}}
\qquad
\text{for all }0\leq j\leq t-1.
\end{equation}
It follows that
\begin{align}
\frac1t\sum_{j=0}^{t-1}\|C_j^t\|_{op}^2
&\leq
\frac{C^2}{t}\sum_{j=0}^{t-1}\frac1{(j+1)^{2(1-\delta)}}\notag\\
&=
\mathcal O\!\left(\frac1{t^{2(1-\delta)}}\right)
=
\mathcal O\!\left(\frac1{t^{1-\delta}}\right),
\label{eq: C-average-square}
\end{align}
where the final, weaker bound is sufficient for our purposes.

Combining
\[
\Phi_j^t=C_j^t+\bar A^{-1}D_j^t
\]
with \eqref{eq: D-average-square} and \eqref{eq: C-average-square}, we
obtain
\begin{equation}\label{eq: Phi-average-square}
\frac1t\sum_{j=0}^{t-1}\|\Phi_j^t\|_{op}^2
=
\mathcal O\!\left(\frac1{t^{1-\delta}}\right).
\end{equation}
Moreover, \eqref{eq: C-bound} and the uniform bound on $D_j^t$ imply
\begin{equation}\label{eq: Phi-uniform}
    \sup_{t\geq1}\max_{0\leq j\leq t-1}\|\Phi_j^t\|_{op}<\infty.
\end{equation}

Finally, we transfer this estimate from $\Phi_j^n$ to $\Upsilon_j^n$.
From the definitions in (\ref{eq: def Upsilon}),
\begin{align*}
\Upsilon_j^n
&=
\epsilon_j\sum_{k=j+1}^{n}
\prod_{l=j+1}^{k-1}(I-\epsilon_l\bar A)
-\bar A^{-1}\\
&=
\frac{\epsilon_j}{\epsilon_{j+1}}
\left(\Phi_{j+1}^{n+1}+\bar A^{-1}\right)
-\bar A^{-1}\\
&=
\frac{\epsilon_j}{\epsilon_{j+1}}\Phi_{j+1}^{n+1}
+
\left(\frac{\epsilon_j}{\epsilon_{j+1}}-1\right)\bar A^{-1}.
\end{align*}
Since
\[
\frac{\epsilon_j}{\epsilon_{j+1}}
=
\left(\frac{j+2}{j+1}\right)^\delta
\]
is uniformly bounded and
\[
\frac{\epsilon_j}{\epsilon_{j+1}}-1
=
\mathcal O\!\left(\frac1{j+1}\right),
\]
\eqref{eq: Phi-uniform}, together with $\Upsilon_n^n=-\bar A^{-1}$, gives
\begin{equation}\label{eq: Upsilon-uniform}
    \sup_{n\geq1}\max_{0\leq j\leq n}\|\Upsilon_j^n\|_{op}<\infty.
\end{equation}
It also follows from \eqref{eq: Phi-average-square} that
\begin{align*}
\frac1n\sum_{j=0}^{n-1}\|\Upsilon_j^n\|_{op}^2
&\leq
\frac{C}{n}\sum_{j=0}^{n-1}\|\Phi_{j+1}^{n+1}\|_{op}^2
+
\frac{C}{n}\sum_{j=0}^{n-1}\frac1{(j+1)^2}\\
&=
\mathcal O\!\left(\frac1{n^{1-\delta}}\right).
\end{align*}
This proves both estimates stated in
\eqref{eq: Upsilon-uniform-and-average}.

\end{document}